\documentclass[a4paper]{cas-sc}
\usepackage[authoryear]{natbib}
\usepackage{mathtools}
\usepackage{hyperref}
\hypersetup{colorlinks, citecolor=teal, linkcolor=Black, urlcolor=Black}
\usepackage{siunitx,booktabs}
\usepackage{amssymb}
\usepackage{siunitx,booktabs}
\usepackage{url}
\usepackage{subcaption}
\usepackage{appendix}
\usepackage{floatrow}
\usepackage{cleveref}
\usepackage{algpseudocode}
\usepackage{algorithm}
\usepackage{amsthm}
\usepackage[normalem]{ulem}

\newtheorem{remark}{Remark}

\def\tsc#1{\csdef{#1}{\textsc{\lowercase{#1}}\xspace}}
\tsc{WGM}
\tsc{QE}
\tsc{EP}
\tsc{PMS}
\tsc{BEC}
\tsc{DE}

\newcommand{\pd}[2]{\frac{\partial #1}{\partial #2}}

\newcommand{\vect}[1]{\mathbf{#1}}
\newcommand{\tens}[1]{\boldsymbol{#1}}

\DeclareSIUnit\bar{bar}
\begin{document}
\shorttitle{Unified Compositional Flow Model for Simulating Multiphase Fractured High Enthalpy System}
\shortauthors{Micheal B. Oguntola, et al.}

\title [mode = title]{Mathematical Modeling of Salt Precipitation and Multi-Phase Flow in High Enthalpy Fractured Geothermal Systems}                      
\author[1]{Micheal B. Oguntola}[orcid = 0000-0001-6692-639X]
\cormark[1]
\ead{micheal.oguntola@uib.no}

\address[1]{University of Bergen, 5007, Bergen, Norway}

\author[1, 2]{Omar Duran}[]
\author[1]{Eirik Keilegavlen}[]
\author[1]{Inga Berre}[]
\address[2]{Stanford University, Stanford, California}

\cortext[cor1]{Corresponding author}
\begin{abstract}
Simulating high-enthalpy fractured geothermal reservoirs is challenging due to the complex coupled processes of non-isothermal, multiphase, multicomponent flow, strongly nonlinear thermodynamics, and  the dominant role of fractures. These complexities are further amplified by mineral scaling, such as halite precipitation, which can severely impair reservoir permeability and well productivity. To address these challenges, this work presents a new compositional flow model, exceeding the capabilities of existing tools. It is based on a persistent set of primary variables (pressure, enthalpy, and overall salt mass fraction) to simulate these coupled processes. The formulation naturally handles phase transitions without manual switching, enhancing numerical stability. The model integrates a discrete fracture-matrix (DFM) approach and employs an efficient and robust correlation-based phase-behaviour linearisation of saltwater thermodynamics, replacing expensive on-the-fly phase separation calculations. It incorporates Kozeny-Carman-type relations to dynamically model porosity and permeability reduction resulting from halite precipitation. Implemented within the open-source PorePy framework, the model is verified through a 1D salt dissolution benchmark against the established closed-source simulator CSMP++, showing strong agreement across geothermal conditions involving transitions between single- and multi-phase regions. Application to a 2D halite-saturated fractured reservoir with injection and production demonstrates the model's capability to predict halite precipitation patterns and their impact on permeability damage and energy recovery. Further, numerical results show the model's value in predicting operational challenges such as wellbore blockage and the role of fracture connectivity. The developed model thus provides an open-source numerical tool for analysing complex heat and mass transport with mineral scaling effects in a high-enthalpy fractured geothermal system. 
\end{abstract}

\maketitle
\section{Introduction}
\label{Introduction}

Meeting the increasing global energy demand while reducing carbon emissions requires expanding renewable energy sources \citep{shukla2022climate, IEA2023}. High-enthalpy geothermal reservoirs can provide clean low-carbon energy with high power generation efficiency \citep{fridleifsson2001geothermal, fridleifsson2008possible,ingridgeothermal,lu2018global}. Such systems are predominantly located along the boundaries of the tectonic plates and are characterised by temperatures that can exceed $200^\circ$C at relatively shallow depths \citep{ingridgeothermal}. At depth, the pore space contains pressurised high-temperature saline fluid consisting primarily of water (H$_2$O) with dissolved salts and non-condensable gases. Sodium-chloride (NaCl) is typically the dominant salt. Based on the pressure-temperature-composition condition, fluid can occur in a supercritical state or as a multiphase mixtures of low-density vapour and higher-density liquid brine, or include precipitated solutes such as halite \citep{afanasyev2018formation}. Continuous energy extraction requires sufficient hydraulic conductivity to facilitate fluid flow between the subsurface heat source and production wells. In low-permeability crystalline and volcanic formations, this conductivity is provided primarily by interconnected fracture networks \citep{tariq2025fracture, lamur2017permeability, zhang2023autonomous, lei2024thermo}, either pre-existing or enhanced through hydraulic stimulation \citep{rahman2002shear, jia2022hydraulic, moska2021hydraulic}. 

During production operations, pressure drawdown near fractures and wellbores can induce boiling of the high-enthalpy geothermal brine \citep{tsypkin2004vapour}. As liquid water vaporises, dissolved salt remains in the residual brine, increasing the local salinity. Conversely, near cold-water injection wells, the inflow of fresh fluid cools the surrounding brine and lowers the local halite solubility. In both cases, once the solubility limit is exceeded, salt precipitates within pore spaces and fractures, impairing the porosity and permeability within the rock formation \citep{phillips1991flow, shahidzadeh2010damage, gunnlaugsson2012scaling, scott2017boiling} and thereby reducing well productivity \citep{hesshaus2013halite, flores2017effect}.

Salt precipitation has caused substantial production losses and in some cases, forced well abandonment \citep{ingridgeothermal} in several major high-enthalpy geothermal fields, including the Salton Sea (USA) \citep{maimoni1982minerals, von2016silica}, Reykjanes (Iceland) \citep{hardardottir2010cu, gunnlaugsson2012scaling, grant2020trace}, and Larderello (Italy) \citep{cavarretta1990schorl}. The effects are particularly severe in fractured reservoirs, where fractures dominate effective permeability despite occupying only a small fraction of the reservoir volume \citep{zhu2021impact, hyman2020flow}. This geometric configuration makes fracture networks susceptible to clogging by precipitated salt, creating flow barriers and altering reservoir connectivity \citep{nooraiepour2018effect}. The impact becomes pronounced in sparse or poorly connected fracture networks, where alternative flow paths between injection and production wells are limited. Furthermore, the heterogeneous nature of fracture networks means that precipitation distributions are highly non-uniform, with certain fractures experiencing preferential clogging depending on local flow rates, temperature gradients, and initial salinity \citep{noiriel2021geometry, chen2024capillary, ji2025capillary}.

To accurately simulate the effect of salt precipitation and dissolution on flow and heat transport in fractured geothermal reservoirs, a suitable fracture model is required to represent the geometric complexity of fractured porous media. Continuum approaches, including dual-porosity and dual-permeability models \citep{gao2022review}, represent fractures through averaged properties distributed over control volumes rather than as explicit geometric entities. These formulations are computationally efficient and widely used in geothermal reservoir simulation \citep{zeng2013numerical, li2019coupled, samardzioska2005numerical, jiang2015multimechanistic}. However, they cannot capture localised precipitation and dissolution, fracture-scale permeability evolution, and the resulting redistribution of flow in reservoirs with highly discontinuous large-scale fractures. In contrast, discrete fracture-matrix (DFM) models explicitly represent fractures as lower-dimensional entities embedded within the surrounding matrix \citep{berre2019flow, martin2005modeling, aghili2021hybrid, xing2017parallel}. This representation preserves fracture connectivity and heterogeneity, making DFM approaches well suited for problems in which localised processes, such as precipitation-induced permeability reduction, strongly impact reservoir behaviour.

Thermal compositional multiphase flow in geothermal reservoirs is governed by strongly coupled mass and energy conservation equations \citep{weis2014hydrothermal, voskov2017operator, chen2006computational,wang2021modeling}. The coupling arises from two main mechanisms. First, phase transition couples the mass and energy equations because variations in phase composition directly affect the thermal state of the fluid mixture. Second, thermodynamic equilibrium relations determine the phase behaviour and fluid properties as functions of primary state variables, such as pressure, enthalpy, and composition \citep{wang2022high}. Consequently, both conservation equations depend on the same thermodynamic description. Closure of the governing equations, therefore, requires constitutive thermodynamic relations and phase-equilibrium conditions, which introduce strong nonlinearities through state-dependent accumulation and flux terms \citep{weis2014hydrothermal}. For binary H$_2$O--NaCl brine, accurately resolving this coupled thermodynamic behaviour is essential to predict phase transitions, heat and salt transport, and precipitation or dissolution processes. 

Solving this strongly coupled, nonlinear system numerically presents significant challenges, particularly when fluid undergoes phase transitions. Traditional compositional simulators commonly employ variable-switching formulations, in which the set of primary unknowns is adjusted locally according to the prevailing fluid phases \citep{alpak2018variable, beaude2019non, aghili2020hybrid, lauser2011new, quiroz2024multi, rajabi2023dynamical}. This strategy can reduce the number of unknowns and improve computational efficiency under stable phase conditions, forming the basis for formulations adopted in widely used geothermal simulators such as TOUGH2 \citep{pruess2003tough}. However, the resulting changes in the algebraic structure introduce non-smoothness that may lead to convergence difficulties near phase boundaries \citep{class2002numerical}. As an alternative, persistent-variable (unified) formulations utilise a fixed set of primary unknowns regardless of phase presence \citep{voskov2012comparison}. By embedding phase behaviour within thermodynamic closure relations and complementarity constraints, these approaches avoid explicit variable switching and provide a more consistent numerical structure across phase transitions, enhancing robustness. The variants of such unified formulations have been implemented for geothermal applications in simulators such as HYDROTHERM \citep{kipp2008guide}, DARTS \citep{khait2019delft, voskov2017operator}, CSMP++ \citep{geiger2006multiphase1,geiger2006multiphase2}, COMPASS \citep{les2025geothermal}, and PorePy \citep{lipovac2025persistent, duran2025mixed, oguntola2025unified}. 

Spatial and temporal discretisation of the governing equations produces a large nonlinear algebraic system that is typically solved using Newton-based methods. In geothermal compositional simulations, a substantial portion of the computational cost arises from repeated evaluation of thermodynamic closure relations and their derivatives during residual and Jacobian assembly. To reduce this computational cost, operator-based linearisation (OBL) strategies have been proposed \citep{voskov2017operator}, in which nonlinear thermodynamic relations are pre-parameterized and represented as interpolated operators in the space of primary state variables. By replacing repeated on-the-fly thermodynamic solves with efficient table lookups and derivative reconstruction, OBL reduces the cost of Jacobian assembly while preserving thermodynamic consistency. Such strategies have been implemented in the DARTS compositional framework \citep{khait2019delft}.

Numerical simulation of brine transport in high-enthalpy geothermal systems requires an accurate characterisation of thermodynamic and transport properties across a wide range of pressures, temperatures, and salinities. In this work, fluid properties for the H$_2$O–NaCl system are determined using the comprehensive correlations developed by \citet{driesner2007system}, which cover magmatic conditions. Their formulations integrate experimental data into a consistent mathematical framework that describes phase stability relations, including the vapour–liquid coexistence surface, halite solubility, and critical curves. One key advantage of these correlations is their ability to efficiently provide direct, non-iterative evaluations of brine fluid properties without the need for computationally intensive  phase separation (flash) calculations and complex equations of state.

Several previous studies have addressed multiphase H$_2$O–NaCl flow in geothermal systems using compositional formulations combined with the H$_2$O–NaCl correlations of \cite{driesner2007system}. \cite{weis2014hydrothermal} developed a framework to simulate thermohaline convection under magmatic-hydrothermal conditions, while \cite{afanasyev2018formation} employed a pressure-enthalpy-salinity formulation to model the formation of brine lenses beneath volcanoes in axisymmetric domains with prescribed high-permeability conduits. \cite{falko2021brine} simulated multiphase thermohaline convection and formation of the brine layer in oceanic crust settings. Although these studies provide useful insights on thermohaline convection and phase separation of H$_2$O-NaCl fluids, they are based mainly on continuum representations or idealised flow conduits and do not address the discrete fracture geometry and connectivity that dominantly
impact salt precipitation and dissolution patterns and the associated permeability evolution in production-driven fractured reservoirs. 

The present work develops a unified multiphase compositional model for simulating brine flow in high-enthalpy fractured geothermal systems, with the central objective of predicting salt precipitation and dissolution patterns and their impact on reservoir performance. The model combines a persistent-variable approach for handling phase appearance and disappearance with a discrete fracture-matrix representation that explicitly resolves the fracture--matrix coupling. Thermodynamic behaviour of the brine system is characterised using the correlations of \citet{driesner2007system}. Following operator-based linearisation principles, thermodynamic properties and their derivatives are precomputed over a discretised primary variable space and retrieved via multilinear interpolation during simulation, eliminating expensive iterative phase separation calculations. The model is implemented in the open-source PorePy framework. Building on recent developments in PorePy, including verification of the unified compositional formulation for non-fractured high-enthalpy systems \citep{oguntola2025unified}, incorporation of discrete fracture-matrix couplings \citep{duran2025mixed}, and implementation of numerical local solvers for analytic flash calculations \citep{lipovac2025persistent}, this work extends the framework to address production-driven flow with salt precipitation and dissolution in high-enthalpy fractured reservoirs. To facilitate reproducibility, we provide not only the numerical model, but also the complete execution workflow in a Docker container \citep{oguntola2026docker}.

The remainder of this paper is organised as follows. Section 2 introduces the mixed-dimensional description of the physical domain using a discrete fracture–matrix (DFM) representation and presents the mathematical formulation, including the governing conservation equations, the persistent primary-variable scheme, and the thermodynamic closure relations. Section 3 describes the numerical discretisation and solution strategy, including space-time discretisation, nonlinear and linear solution methods, and secondary variable evaluation scheme. Section 4 presents verification against the established simulator CSMP++ \citep{weis2014hydrothermal} for a 1D problem involving halite (solid salt) dissolution and multiphase behaviour; additional verification studies under high-enthalpy conditions are detailed in \citet{oguntola2025unified}. Section 5 applies the framework to 2D fractured reservoir cases, examining how fracture connectivity and aperture–halite feedback together control precipitation and dissolution patterns and their impacts on reservoir productivity. Section 6 concludes with a summary of key findings.

\section{Mathematical model}
This section presents the mathematical framework for simulating non-isothermal multiphase compositional flow in high-enthalpy fractured H$_2$O-NaCl systems with salt precipitation and dissolution effects. We introduce a mixed-dimensional representation of the physical domain in which fractures and their intersections are treated as lower-dimensional manifolds embedded in the porous matrix, and formulate conservation equations in fractional form for mass and energy using a persistent set of primary variables--pressure, enthalpy, and overall salt mass fraction--that remain fixed across all phase transitions. The section concludes with constitutive relations that close the equation system and a specification of initial and boundary conditions required for a well-posed formulation.

\subsection{Mixed-dimensional geometry}
\label{subsec:mixed_dim_geometry}
\begin{figure}[htp!]
    \centering
    \includegraphics[width=0.8\linewidth]{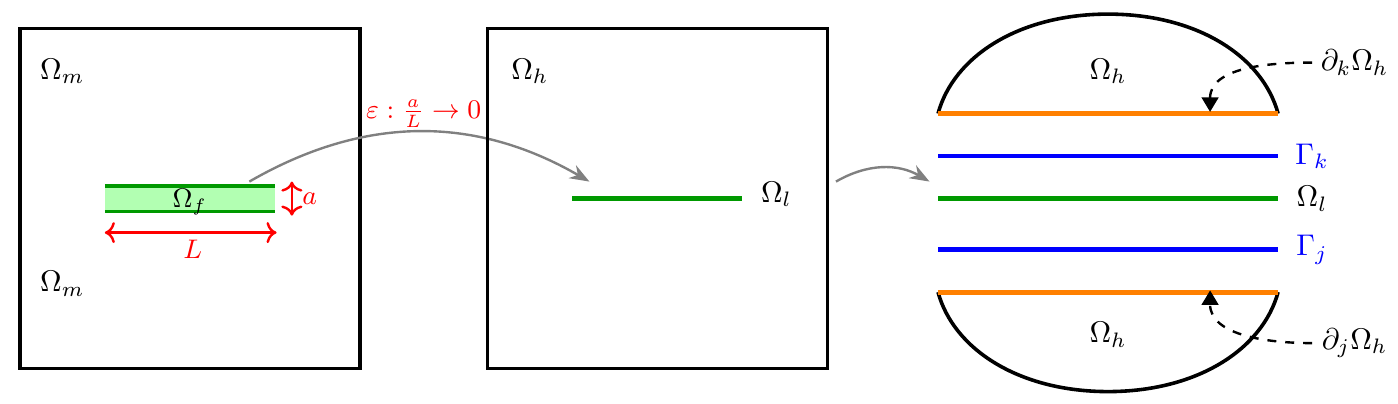}
    \caption{
        Dimensional reduction of a 2D fractured porous medium. Left: Equidimensional representation, where a thin fracture domain $\Omega_f\subset\mathbb{R}^2$ of aperture $a$ and characteristic length $L$ is embedded in a porous matrix $\Omega_m\subset\mathbb{R}^2$. Middle: Mixed-dimensional representation, where the fracture aspect ratio $\varepsilon$ is considered very small. Consequently, the 2D fracture is reduced to a 1D manifold $\Omega_l$ (lower-dimensional subdomain) embedded within the 2D matrix $\Omega_h$ (higher-dimensional subdomain). Right: Illustration showing the reduced fracture $\Omega_l$, the adjacent matrix subdomains $\Omega_h$, and the corresponding matrix--fracture interfaces $\Gamma_{j}$ and $\Gamma_{k}$, together with the internal matrix boundaries $\partial_{j}\Omega_h$ and $\partial_{k}\Omega_h$. The subscripts $j$ and $k$ denote interfaces on the two sides of the fracture, which geometrically coincide with $\partial_{j}\Omega_h$ and $\partial_{k}\Omega_h$ respectively. Adapted from \citet{keilegavlen2021porepy}.
    }
    \label{fig:mat-frac-mdg}
\end{figure}

\begin{figure}[htp!]
    \centering
    \includegraphics[width=0.8\linewidth]{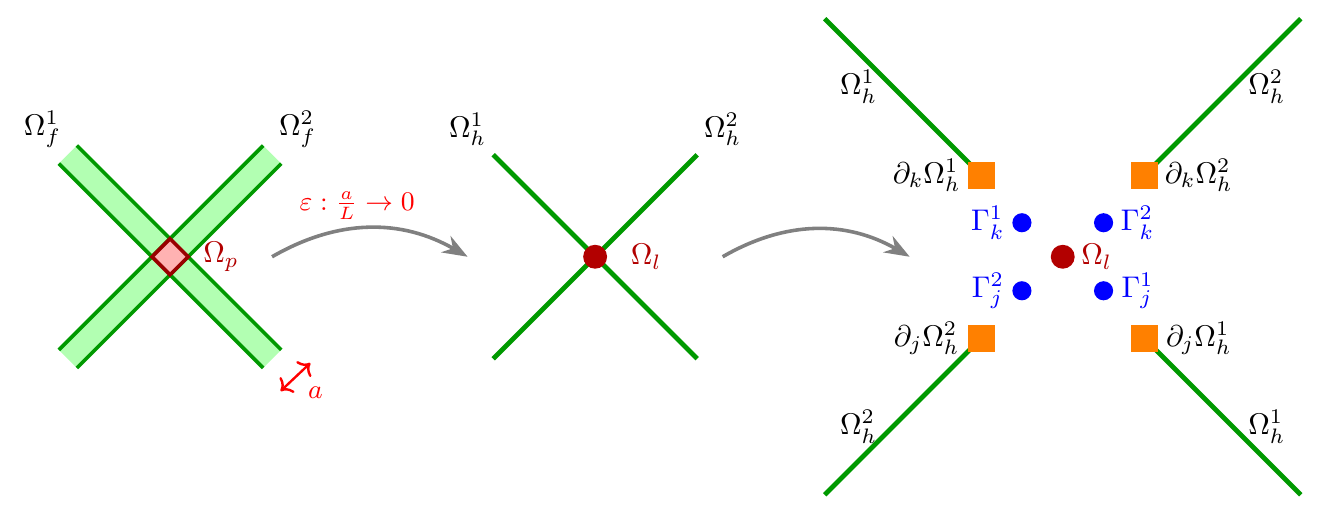}
    \caption{Dimensional reduction of intersecting 2D fractures.
        Left: Equidimensional configuration where two fracture domains $\Omega_f^i\subset \mathbb{R}^2, i=1,2$ of aperture $a$ and characteristic length $L$ intersect, forming an intersection region $\Omega_p\subset\mathbb{R}^2$. Middle: Mixed-dimensional representation, considering a very small fracture aspect ratio $\varepsilon,$ where the 2D fractures are reduced to 1D manifolds $\Omega_h^i, i=1,2$ and the intersection region is reduced to a 0D point $\Omega_l.$ Here, the reduced 1D fractures are the higher-dimensional subdomain relative to the intersection point, which acts as the lower-dimensional subdomain. Right: Illustration showing the fracture subdomains $\Omega_h^i, i=1,2$, the fracture internal boundaries $\partial_{j}\Omega_h^i, \partial_{k}\Omega_h^i$, and the fracture--point interfaces $\Gamma_{j}^i, \Gamma_{k}^i$. The subscripts $j$ and $k$ denote the interfaces where the manifolds geometrically coincide with the fracture internal boundaries. Adapted from \citet{keilegavlen2021porepy}.
    }
    \label{fig:frac-iter-mdg}
\end{figure}

We consider a 2D domain $\Omega \subset \mathbb{R}^2$ composed of a porous rock matrix containing a network of fractures. Traditionally, in an equi-dimensional representation, each fracture $\Omega_f$ and fracture intersection $\Omega_p$ have the same dimension as the surrounding matrix (
Figures~\ref{fig:mat-frac-mdg} and~\ref{fig:frac-iter-mdg}, left). However, when the aperture-to-length (or aspect) ratio $\varepsilon$ is sufficiently small, explicitly resolving the full fracture geometry is computationally expensive, and variable variation across the fracture width is often negligible \citep{dugstad2022dimensional}. This motivates a dimensional reduction in which fractures and their intersections are represented as lower-dimensional manifolds (Figures~\ref{fig:mat-frac-mdg} and~\ref{fig:frac-iter-mdg}, middle). 

Following the discrete fracture-matrix (DFM) approach \citep{martin2005modeling, boon2018robust, berre2019flow, keilegavlen2021porepy}, we decompose $\Omega$ into subdomains of different dimensions. Let $\Omega_d$ denote the collection of all subdomains of dimension $d$, so that

\begin{equation}
    \Omega = \bigcup_{d=0}^{2} \Omega_d, \qquad
    \Omega_d = \bigcup_{i\in\mathcal{I}_d} \Omega_d^i,
    \label{eq:domain_decomposition}
\end{equation}
where $\Omega_2$ represents the 2D matrix, $\Omega_1$ the set of 1D fracture subdomains, $\Omega_0$ the set of 0D fracture intersection subdomains, and $\mathcal{I}_d$ the index set identifying the individual subdomains of dimension $d$. The coupling between subdomains of adjacent dimensions is obtained through interfaces. For a fracture $\Omega_l$ embedded in a matrix $\Omega_h$, we define the upper and
lower interfaces $\Gamma_k$ and $\Gamma_j$, each inheriting the geometry of the fracture and coinciding geometrically with the matrix internal boundaries $\partial_k\Omega_h$ and $\partial_j\Omega_h$ (Figure~\ref{fig:mat-frac-mdg}, right). 

\begin{figure}[htp!]
    \centering
    \includegraphics[width=0.5\linewidth]{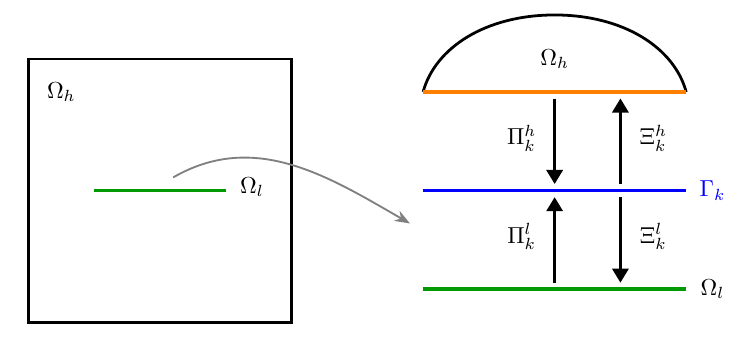}
   \caption{Matrix--fracture coupling structure. Left:
Mixed-dimensional geometry where the fracture subdomain $\Omega_l$ is embedded in the matrix subdomain $\Omega_h.$ Right: The interface $\Gamma_k$ between $\Omega_h$ and $\Omega_l$, with trace operators $\Pi_k^h$, $\Pi_k^l$ mapping subdomain quantities to the interface  and projection operators $\Xi_k^h$, $\Xi_k^l$ mapping interface quantities back to the subdomains. The coupling on the opposite side is analogous. Adapted from \citet{keilegavlen2021porepy}.}
\label{fig:mat-frac-mdg-c} 
\end{figure}

The mapping of quantities between subdomains and interfaces is handled by specific operators (see Figure \ref{fig:mat-frac-mdg-c}): the trace operators $\Pi_k^h$ and $\Pi_k^l$ map subdomain quantities to the interfaces, while the prolongation (or projection) operators $\Xi_k^h$ and $\Xi_k^l$ map interface quantities back to the adjacent subdomains \citep{keilegavlen2021porepy}. With these operators, we are able to directly quantify fluxes across subdomains of codimension one. The geometric reduction from the equi-dimensional fractures and intersections to lower-dimensional manifolds introduces a specific volume $\nu_d$ that accounts for the average across the reduced dimensions:

\begin{equation}
    \nu_d := a^{n-d} = 
    \begin{cases}
        1   & \text{for } d = 2 \quad \text{(matrix)}, \\[2pt]
        a   & \text{for } d = 1 \quad \text{(fracture)},\\[2pt]
        a^2  & \text{for } d = 0 \quad \text{(fracture intersection)},
    \end{cases}
    \label{eq:specific_volume_term}
\end{equation}
where $n=2$ is the ambient dimension, $d$ is the subdomain dimension, and $a$ is the local effective aperture. For an intersection point, $a$ is the average aperture of the intersecting fractures \citep{stefansson2021fully}. This factor ensures dimensional consistency when conservation laws, originally formulated over the equi-dimensional domains, are expressed on the reduced geometry. The framework extends naturally to 3D, where the dimensional hierarchy spans $d = 0, \ldots, 3$ \citep{boon2018robust}.


\subsection{Model equations}
\label{sec:model_equations}
Building on the mixed-dimensional representation of the physical domain $\Omega$ introduced in the previous subsection, we extend the thermal multiphase compositional model, in fractional form, for a continuum domain presented in \citep{oguntola2025unified} to the mixed-dimensional setting. This extension introduces the specific volume
$\nu_d$ to account for the dimensional reduction, together with interface
coupling terms that enforce conservation of mass and energy across
subdomains of codimension one. Under no gravity and negligible capillarity conditions, we assume that the reservoir fluid consists of
$N_c = 2$ chemical components--water (H$_2$O) and salt (NaCl), indexed
by $\xi$--which can coexist in up to $N_p = 3$ phases--liquid brine (liq), vapour (vap), and solid halite (hal), indexed
by $\gamma$. Precipitated halite is treated as an immobile phase that contributes to mass and energy storage but not to flow, hence, its relative permeability is set to zero. For each subdomain $\Omega_d^i$, $i \in \mathcal{I}_d$, of dimension $d \in \{0, 1, 2\}$ and specific volume $\nu_d$ given by Equation~\eqref{eq:specific_volume_term}, we obtain the mixed-dimensional balance equations by integrating their
equi-dimensional counterparts \citep{oguntola2025unified} across the reduced dimension, under the assumption that solution variables vary negligibly in the direction normal to the fracture \citep{dugstad2022dimensional, martin2005modeling}. These equations enforce conservation of mass for each component and conservation of energy for the fluid--rock system. 

Throughout this section, the index $i$ identifies variables defined on subdomain $\Omega_d^i$; definitions and units of all quantities are listed in Tables~\ref{tab:primary_variables}--\ref{tab:other_physical_quantities}, where in most cases the index $i$ is suppressed for readability.

\subsubsection{Component mass balance}
\label{component-mass-balance}
For each component $\xi \in \mathrm{\{ H_2 O, NaCl\}}$ the overall mass balance equation on a subdomain $\Omega_d^i$ is given by
\begin{equation}
	\pd{}{t}\bigl(\nu_d \phi_i\rho_i z_i^\xi \bigr)
	+ \nabla \cdot \bigl(\nu_d \vect{F}_i^\xi \bigr)
    - \sum_{k\in\hat{\mathcal{S}}_i} \Xi_k^i \nu_kF_k^\xi
	= \nu_d \, q_i^\xi,
	\label{eq:component_mass_balance}
\end{equation}
where $\phi_i$ is the porosity, $\rho_i = \sum_\gamma s_i^\gamma \rho_i^\gamma$ is the fluid mixture density, in which $s_i^\gamma$ and $\rho_i^\gamma$ are the saturation and density of phase $\gamma,$ and $z_i^\xi = \frac{1}{\rho_i}\sum_\gamma s_i^\gamma\rho_i^\gamma \chi_{\xi\gamma}^i$ is the overall mass fraction of component $\xi,$ in which $\chi_{\xi\gamma}^i$ is the partial mass fraction of component $\xi$ in phase $\gamma.$ The projection operator $\Xi_k^i$ maps quantities from interface $\Gamma_k$ to the subdomain $\Omega_d^i,$ the interface specific volume is $\nu_k = \Pi_k^h \nu_{h},$ and $\hat{\mathcal{S}}_i$ is the set of all interfaces connecting $\Omega_d^i$ to higher-dimensional subdomains. $q_i^\xi$ is a volumetric source or sink term defined per $n$-dimensional volume.

The component mass flux $\vect{F}_i^\xi$, in fractional form, is given by \citep{chen2006computational, oguntola2025unified}
\begin{equation}
    \vect{F}_i^\xi = - f_i^\xi\lambda_i\tens{K}_i\nabla p_i,
    \label{eq:component_flux_diffusive}
\end{equation}
where $p_i$ is the fluid pressure and $\tens{K}_i$ is the rock permeability tensor. The component fractional flow is $f_i^\xi = \frac{\lambda_i^\xi}{\lambda_i}$, with the component density-weighted mobility $\lambda_i^\xi = \sum_{\gamma} \lambda_{\xi\gamma}^i,$ where $\lambda_{\xi\gamma}^i = \frac{\rho_i^\gamma\chi_{\xi\gamma}^i k_{r\gamma}^i}{\mu_i^\gamma}$ is the component mobility in a phase $\gamma,$ with $k_{r\gamma}^i$ and $\mu_i^\gamma$ as the relative permeability and viscosity of phase $\gamma,$ respectively. The total density-weighted mobility is $\lambda_i = \sum_\xi \lambda_i^\xi.$ Here, we note that only mobile phases contribute to the flux $\vect{F}_i^\xi,$ as precipitated halite is immobile.

Equivalently, Equation~\eqref{eq:component_flux_diffusive} can be written as
\begin{equation}
	\vect{F}_i^\xi = f_i^\xi \, \lambda_i \, \vect{u}_i,
	\label{eq:component_flux_advective}
\end{equation}
where 
\begin{equation}
    \vect{u}_i = - \tens{K}_i \nabla p_i,
    \label{eq:darcy-flux}
\end{equation}
is referred to in this work as mobility-free Darcy flux (i.e. the usual single-phase Darcy flux \citep{chen2006computational}).
\begin{remark}
    The splitting of $\vect{F}_i^\xi$ into $\vect{u}_i$ and the mobility term $f_i^\xi \lambda_i$ in Equation~\eqref{eq:component_flux_advective} is motivated by the numerical discretisation of the advective term in Equation \eqref{eq:component_mass_balance} using an upwind scheme. In this case, $\vect{u}_i$ is discretised independently (e.g., by multi-point flux approximation), while $f_i^\xi \lambda_i$ is evaluated based on the direction of $\vect{u}_i$. An equivalent formulation instead incorporates the total mobility $\lambda_i$ into the diffusive tensor of the Darcy flux, in which case only the dimensionless fractional weight $f_i^\xi$ is upwound. Both discretisations are consistent approximations of the same continuous flux \eqref{eq:component_flux_advective}, differing in the discrete treatment of the mobility.
\end{remark}

The interface component mass flux is
\begin{equation}
	F_k^\xi = f_k^\xi \, \lambda_k \, u_k,
	\label{eq:interface_component_flux}
\end{equation}
where $u_k$ is herein referred to as the interface mobility-free Darcy flux given by \citet{martin2005modeling}
\begin{equation}
    u_k = -K_k \, \frac{2}{\Pi_k^l a_l}
      \bigl(\Pi_k^l p_l - \Pi_k^h p_h\bigr),
	\label{eq:interface_darcy}
\end{equation}
in which $K_k$ is the scalar normal permeability on $\Gamma_k$ (inherited
from the lower-dimensional subdomain), $p_l$ and $p_h$ are the lower- and
higher-dimensional subdomain pressures, and $a_l$ is the aperture of the
lower-dimensional subdomain. The quantities $f_k^\xi$ and $\lambda_k$ denote
the interface component fractional flow and total mobility, respectively.


\subsubsection{Total mass balance}
\label{subsec:total_mass_balance}
Summing Equation~\eqref{eq:component_mass_balance} over all components $\xi$
and using $\sum_\xi f_i^\xi = 1$, $\sum_\xi z_i^\xi = 1,$ and Equation \eqref{eq:component_flux_advective} yields the total
mass balance (or pressure equation) on the subdomain $\Omega_d^i$:
\begin{equation}
    \pd{}{t}\bigl(\nu_d \, \phi_i \rho_i\bigr)
    + \nabla \cdot \bigl(\nu_d \, \vect{F}_i\bigr)
    - \sum_{k \in \hat{\mathcal{S}}_i} \Xi_k^i \, \nu_k \, F_k
    = \nu_d \, q_i,
    \label{eq:total_mass_balance}
\end{equation}
where $\vect{F}_i = \lambda_i \, \vect{u}_i$ is the total mass flux, $q_i = \sum_\xi q_i^\xi$ is the total source or sink, $F_k = \sum_\xi F_k^\xi$ is the interface total mass flux, and $\vect{u}_i$ is given by Equation \eqref{eq:darcy-flux}.

\subsubsection{Energy balance}
\label{subsec:energy_balance}
Using the assumption of local thermal equilibrium \citep{oguntola2025unified, chen2006computational}, the energy balance on a subdomain $\Omega_d^i$ in terms of the specific enthalpy of fluid mixture is given as:
\begin{equation}
    \pd{}{t}\bigl(\nu_d \, E_i\bigr)
    + \nabla \cdot \bigl[\nu_d
      \bigl(\vect{w}_i + \vect{q}_i\bigr)\bigr]
    - \sum_{k \in \hat{\mathcal{S}}_i} \Xi_k^i \, \nu_k
      \bigl(w_k + q_k\bigr)
    = \nu_d \, q_i^E,
    \label{eq:energy_balance}
\end{equation}
where $q_i^E$ is an energy source or sink term. The energy accumulation term includes contributions from all fluid phases (including precipitated halite) and the rock matrix:
\begin{equation}
    E_i = \phi_i (\rho_i h_i - p_i) + (1-\phi_i)\rho_r c_{p,r} T_i,
    \label{eq:energy_accumulation}
\end{equation}
where $h_i = \frac{1}{\rho_i}\sum_{\gamma} s_i^\gamma\rho_i^\gamma h_i^\gamma$ is the specific enthalpy of fluid mixture, $\rho_r$ is the rock density, $c_{p,r}$ is the rock specific heat capacity, and $T_i$ is the temperature. As with the mass accumulation, the summation over all
phases ensures that the thermal energy stored in precipitated halite is
accounted for. The subdomain energy flux comprises an advective enthalpy flux $\vect{w}_i$  and a conductive heat flux $\vect{q}_i$:
\begin{equation}
    \vect{w}_i = \sum_{\gamma}
      \lambda_i^\gamma h_i^\gamma \vect{u}_i, \qquad
    \vect{q}_i = -\tens{D}_{i}^h \nabla T_i,
    \label{eq:energy_fluxes}
\end{equation}
where $h_i^\gamma$ is the specific enthalpy of phase $\gamma$, $\lambda_i^\gamma = \sum_\xi \lambda_{\xi\gamma}^i$ is the phase density-weighted mobility, and the effective thermal conductivity tensor is $\tens{D}_{i}^h = \bigl(\phi_i \sum_{\gamma}
s_i^\gamma \kappa^\gamma + (1 - \phi_i) \kappa_r\bigr) \tens{I}_n$, with $\kappa^\gamma$ and $\kappa_r$ the thermal conductivities of phase $\gamma$ and the solid rock, respectively. The mobility-free Darcy flux $\vect{u}_i$ is given by Equation \eqref{eq:darcy-flux}. As with the mass flux in Equation \eqref{eq:component_flux_advective}, only mobile phases contribute to the advective enthalpy flux.

The interface energy flux on $\Gamma_k$ comprises an advective contribution $w_k$ and a conductive contribution $q_k$:
\begin{equation}
    w_k = \sum_{\gamma}
      \lambda_k^\gamma h_k^\gamma \, u_k, \qquad
    q_k = -D_k^h \, \frac{2}{\Pi_k^l a_l}
      \bigl(\Pi_k^l T_l - \Pi_k^h T_h\bigr),
    \label{eq:interface_energy_fluxes}
\end{equation}
where $u_k$ is the interface mobility-free Darcy flux defined in Equation~\eqref{eq:interface_darcy}, $D_k^h$ is the interface thermal conductivity, $T_l$ and $T_h$ are the lower- and higher-dimensional subdomain temperatures, and $\lambda_k^\gamma$ and $h_k^\gamma$ are the interface density-weighted mobility and specific enthalpy of phase $\gamma.$

\begin{remark}\leavevmode
\begin{enumerate}
    \item In Equations~\eqref{eq:component_mass_balance},
    \eqref{eq:total_mass_balance}, and~\eqref{eq:energy_balance}: for
    $d = 0$ (intersection subdomains), the divergence term vanishes since a point has no spatial extent; for $d = 1$ (fracture subdomains), the divergence operator reduces to its tangential component along the
    fracture; and for $d = 2$ (matrix subdomain), the interface source term vanishes since the matrix has no higher-dimensional neighbour.
    \label{rm:subdomain_interface_equations}
    \item Mass and energy fluxes between $\Omega_d^i$ and its lower-dimensional neighbours of codimension one are imposed as internal boundary conditions on $\Omega_d^i$, enforcing continuity of normal fluxes at the internal boundaries $\partial_\alpha \Omega_d^i$ that coincide geometrically with the interfaces $\Gamma_\alpha \in \check{\mathcal{S}}_i.$ Here, $ \check{\mathcal{S}}_i$ is the set of interfaces connecting $\Omega_d^i$ to the lower-subdomains.
    \label{rm:internal_flux_continuity}
\end{enumerate}
\end{remark}

\subsection{Persistent-variable formulation}
Equations \eqref{eq:component_mass_balance}, \eqref{eq:total_mass_balance}, and \eqref{eq:energy_balance} form a strongly coupled system of nonlinear mixed-dimensional PDEs. A critical challenge in using this system to simulate high-enthalpy fractured geothermal reservoirs is the robust treatment of phase transitions. 
To this end, we adopt the persistent-variable formulation \citep{oguntola2025unified, duran2025mixed, lipovac2025persistent}, in
which a fixed set of primary variables is maintained regardless of the number or type of fluid phases present.

\subsubsection*{Choice of primary variables}
The primary variables are selected to span the thermodynamic state space continuously across all phase configurations:
\begin{equation}
    \vect{x}_i = (p_i, \, h_i, \, z_i),
    \label{var:primary_variables}
\end{equation}
where $p_i$ is the fluid pressure, $h_i$ is the specific enthalpy of the fluid
mixture, and $z_i = z_i^{\text{NaCl}}$ is the overall mass fraction of salt on each subdomain $\Omega_d^i$. Since, constitutively, $\sum_\xi z^\xi_i = 1$, the water mass fraction is determined by $z^{\mathrm{H_2O}}_i = 1 - z_i$. The pressure $p_i$ is solved for by the total mass balance \eqref{eq:total_mass_balance}, the salt fraction $z_i$ by the component mass balance \eqref{eq:component_mass_balance}, and the enthalpy $h_i$ by the energy balance \eqref{eq:energy_balance}.
\begin{table}[hpt!]
    \centering
    \caption{Definition of primary variables.}
    \label{tab:primary_variables}
    \begin{tabular}{lll}
        \toprule
        \textbf{Variable} & \textbf{Meaning} & \textbf{Units} \\
        \midrule
        $p$   & Fluid pressure                        & Pa \\
        $h$   & Specific enthalpy of fluid mixture     & J/kg \\
        $z$   & Overall mass fraction of salt (NaCl)   & -- \\
        \bottomrule
    \end{tabular}
\end{table}

\subsection{Constitutive relations}
\label{constitutive_relations}
The balance equations
\eqref{eq:component_mass_balance}--\eqref{eq:energy_balance} contain
secondary variables $\vect{y}_i$ given by Table~\ref{tab:secondary_variables}
and state-dependent material properties such as porosity and permeability
(see Table~\ref{tab:other_physical_quantities}) that must be eliminated by functions of the
primary variables to close the system. Given the primary variables \eqref{var:primary_variables}, the secondary variables are recovered
through the map
\begin{equation}
    \vect{y}_i = \vect{y}(\vect{x}_i),
    \label{eq:secondary_recovery}
\end{equation}
where $\vect{y}_i = (T_i, s_i^\gamma, \chi_{\xi\gamma}^i, \rho_i^\gamma,
\mu_i^\gamma, h_i^\gamma, k_{r\gamma}^i)$. For the H$_2$O--NaCl system, the thermodynamic variables---temperature, phase saturations, compositions, densities, viscosities, and enthalpies---are evaluated using the saltwater correlation formulae of
\citet{driesner2007system}, valid for pressures up to 5000~bar, temperatures up to 1000$^\circ$C, and salt mass fractions across the full range $z \in [0, 1]$. Since these formulae are usually given in the $(p, T, z)$ space, we consider their reconstruction in the primary variables \eqref{var:primary_variables} space using a bisection method to ensure thermodynamic consistency. Details of their efficient implementation are given in \citet{oguntola2025unified}.

When halite precipitates as an immobile solid phase, it occupies a fraction of the pore space. The effective porosity of the matrix subdomain $(d = 2)$ available to mobile fluid
phases becomes \citep{weis2014hydrothermal}
\begin{equation}
    \phi_i = \phi^0 - s_i^{\text{hal}} \phi^0 
           = \phi^0 (1 - s_i^{\text{hal}}),
    \label{eq:porosity_update}
\end{equation}
and the intrinsic permeability is updated accordingly through the Kozeny-Carman relation:
\begin{equation}
    \tens{K}_i = \tens{K}^0 \left(\frac{\phi_i}{\phi^0}\right)^2
    \label{eq:permeability_kozeny_carman}
\end{equation}
where $\tens{K}^0$ is the reference permeability corresponding to the
reference porosity $\phi^0$. In fracture and intersection subdomains, the aperture also evolves with halite precipitation according to
\begin{equation}
    a_i = \max\Bigl(a^0 \bigl(1 - s_i^{\text{hal}}\bigr)^\varphi, \; a_{\min}\Bigr),
    \label{eq:aperture_update}
\end{equation}
where $a^0$ is the reference fracture aperture, $a_{\min}$ is a minimum aperture imposed to avoid vanishing permeability, and $\varphi\geq 0$ is an exponent to control the clogging scale by precipitated halite. The fracture permeability is then given by the cubic law:
\begin{equation}
    \mathcal{K}_i = \frac{a_i^2}{12} \textbf{I}_n.
    \label{eq:cubic_law}
\end{equation}

The relative permeability of each mobile phase follows a modified Corey-type
model that accounts for pore-space reduction by halite precipitation, similar to the formulation in \citet{falko2021brine, weis2014hydrothermal}.

\subsection{Initial and boundary conditions}
\label{sec:ic_bc}

To complete the model formulation, we specify initial and boundary conditions 
for the coupled system of balance equations 
\eqref{eq:component_mass_balance}--\eqref{eq:energy_balance} on the mixed-dimensional 
domain $\Omega$.

\subsubsection*{Initial conditions}
At time $t = 0$, the thermodynamic state of each subdomain $\Omega_d^i$ is prescribed 
through the primary variables:
\begin{equation}
    \mathbf{x}_i\big|_{t=0} 
    = \mathbf{x}_{i,0} 
    := \bigl(p_{i,0},\; h_{i,0},\; z_{i,0}\bigr) 
    \qquad \text{on } \Omega_d^i,
    \label{eq:ic}
\end{equation}
where $p_{i,0}$, $h_{i,0}$, and $z_{i,0}$ denote the initial distributions of pressure, specific enthalpy of the fluid mixture, and overall salt mass fraction, respectively. All initial secondary variables are determined from $\mathbf{x}_{i,0}$ through the thermodynamic closure map~\eqref{eq:secondary_recovery}. When precipitated halite is present at the initial time (i.e.,$s_i^{\mathrm{hal}}\big|_{t=0} > 0$), the initial porosity, permeability, and fracture aperture are computed from Equations~\eqref{eq:porosity_update}--\eqref{eq:cubic_law} accordingly.

\subsubsection*{Boundary conditions}
Let $\partial \Omega_{\mathrm{ext}}^i$ and $\partial \Omega_{\mathrm{int}}^i$ denote the external and internal boundaries of the subdomains $\Omega_d^i, d=1, 2$. We decompose $\partial \Omega_{\mathrm{ext}}^i$ into non-overlapping Dirichlet $\partial\Omega_{\mathrm{ext}}^{i, D}$ and Neumann $\partial\Omega_{\mathrm{ext}}^{i, N}$ portions.

On Dirichlet boundaries, the full set of primary variables is prescribed:
\begin{equation}
	p_i = p_D, \qquad h_i = h_D, \qquad z_i = z_D 
	\qquad \text{on } \partial \Omega_{\mathrm{ext}}^{i, D},
	\label{eq:bc_dir}
\end{equation}
which fully determines the thermodynamic state at the boundary through 
the closure map~\eqref{eq:secondary_recovery}. 

On Neumann boundaries, the normal component mass fluxes and energy flux are prescribed:
\begin{equation}
	\nu_d\, \mathbf{F}_i^\xi \cdot \mathbf{n} = Q_N^\xi, 
	\quad \xi \in \{\text{H}_2\text{O},\, \text{NaCl}\},
	\qquad
	\nu_d\, \bigl(\mathbf{w}_i + \mathbf{q}_i\bigr) 
	\cdot \mathbf{n} = q_N^{E}
	\qquad \text{on } \partial \Omega_{\mathrm{ext}}^{i,\, N},
	\label{eq:bc_neu}
\end{equation}
where $\mathbf{n}$ is the outward unit normal to $\partial\Omega_{\mathrm{ext}}^{i,\, N}$, $Q_N^\xi$ is a prescribed normal component mass flux and $q_N^{E}$ is a prescribed normal energy flux. The total mass flux on $\partial \Omega_{\mathrm{ext}}^{i,\, N}$ is given consistently by $Q_N = \sum_\xi Q_N^\xi.$ No-flow and adiabatic conditions correspond to $Q_N^\xi = 0$ for all $\xi$ and $q_N^{E} = 0$, respectively.

On the internal boundaries $\partial \Omega_{\mathrm{int}}^i 
= \bigcup_{\alpha \in \check{\mathcal{S}}_i} 
\partial_\alpha \Omega_d^i$, which coincide geometrically with the interfaces $\Gamma_\alpha$ connecting $\Omega_d^i$ to its lower-dimensional neighbours 
(cf.\ Remark~\ref{rm:internal_flux_continuity}), conservation of mass and 
energy is enforced by requiring continuity of normal fluxes:
\begin{equation}
	\nu_d\, \vect{F}_i \cdot \mathbf{n}_\alpha 
	=  \Xi_\alpha^i\,\nu_\alpha\, F_\alpha, \qquad
	\nu_d\, \bigl(\mathbf{w}_i + \mathbf{q}_i\bigr) 
	\cdot \mathbf{n}_\alpha 
	=  \Xi_\alpha^i\,\nu_\alpha\, \bigl(w_\alpha + q_\alpha\bigr)
	\qquad \text{on } \partial_\alpha \Omega_d^i, 
	\quad \alpha \in \check{\mathcal{S}}_i,
	\label{eq:bc_int}
\end{equation}
where $\mathbf{n}_\alpha$ is the outward unit normal to $\partial_\alpha \Omega_d^i$, $\Xi_\alpha^i$ is the projection operator mapping quantities from interface $\Gamma_\alpha$ to the subdomain $\Omega_d^i$, $\nu_\alpha$ is the interface specific volume, $F_\alpha$ is the interface total mass flux, and $w_\alpha, q_\alpha$ are the interface advective enthalpy and conductive heat fluxes.
\begin{remark}
	Lower-dimensional subdomains $\Omega_d, d = 0, 1$ that lie in the interior of $\Omega_d, d=2$ receive no external boundary conditions; their coupling to adjacent subdomains is governed entirely by the interface flux terms in Equations~\eqref{eq:component_mass_balance}--\eqref{eq:energy_balance}.
\end{remark}

\begin{table}[H]
    \centering
    \caption{Definition of secondary variables.}
    \label{tab:secondary_variables}
    \begin{tabular}{lll}
        \toprule
        \textbf{Variable} & \textbf{Meaning} & \textbf{Units} \\
        \midrule
        $T$                & Temperature                                         & K \\
        $s^\gamma$         & Volumetric saturation of phase $\gamma$             & -- \\
        $\chi_{\xi\gamma}$ & Mass fraction of component $\xi$ in phase $\gamma$  & -- \\
        $\rho^\gamma$      & Density of phase $\gamma$                           & kg/m$^3$ \\
        $\mu^\gamma$       & Dynamic viscosity of phase $\gamma$                 & Pa$\cdot$s \\
        $h^\gamma$         & Specific enthalpy of phase $\gamma$                 & J/kg \\
        \bottomrule
    \end{tabular}
\end{table}
\begin{table}
    \centering
    \caption{Definition of other physical quantities.}
    \label{tab:other_physical_quantities}
    \begin{tabular}{>{\raggedright}p{2.2cm}p{7.5cm}l}
        \toprule
        \textbf{Symbol} & \textbf{Meaning} & \textbf{Units} \\
        \midrule
        \multicolumn{3}{l}{\textit{Derived fluid quantities}} \\
        $\rho$              & Fluid mixture density
        & kg/m$^3$ \\
        $k_{r\gamma}$       & Relative permeability of phase $\gamma$             & -- \\
        $\lambda_{\xi\gamma}$ & Density-weighted mobility of $\xi$ in $\gamma$ & Pa$^{-1}\!\cdot$s \\
        $\lambda^\xi$       & Component density-weighted mobility & Pa$^{-1}\!\cdot$s \\
        $\lambda$           & Total density-weighted mobility & Pa$^{-1}\!\cdot$s \\
        $\lambda^\gamma$    & Phase density-weighted mobility & Pa$^{-1}\!\cdot$s \\
        $f^\xi$             & Component fractional flow
        & -- \\
        \midrule
        \multicolumn{3}{l}{\textit{Rock and thermal properties}} \\
        $\phi$              & Porosity                     & -- \\
        $\phi^0$            & Reference porosity (prior to precipitation)  & -- \\
        $\tens{K}$          & Intrinsic permeability tensor & m$^2$ \\
         $\mathcal{K}$ & Cubic law fracture permeability & m$^2$ \\
        $\tens{K}^0$          & Reference permeability tensor & m$^2$ \\
        $K_k$          & Scalar normal permeability on $\Gamma_k$
        & m$^2$ \\
        $\rho_r$            & Rock density                 & kg/m$^3$ \\
        $c_{p,r}$           & Rock specific heat capacity  & J$\cdot$kg$^{-1}\!\cdot$K$^{-1}$ \\
        $\kappa_r$          & Rock thermal conductivity    & W$\cdot$m$^{-1}\!\cdot$K$^{-1}$ \\
        $\kappa^\gamma$     & Thermal conductivity of phase $\gamma$
        & W$\cdot$m$^{-1}\!\cdot$K$^{-1}$ \\
        $D_k^h$        & Interface thermal conductivity
        & W$\cdot$m$^{-1}\!\cdot$K$^{-1}$ \\
        $\tens{D}^h$        & Effective thermal conductivity tensor
        & W$\cdot$m$^{-1}\!\cdot$K$^{-1}$ \\
        \midrule
        \multicolumn{3}{l}{\textit{Fluxes}} \\
        $\vect{F}^\xi$      & Total component mass flux (subdomain)
        & kg$\cdot$m$^{-2}\!\cdot$s$^{-1}$ \\
        $F_k^\xi$           & Interface component mass flux
        & kg$\cdot$m$^{-2}\!\cdot$s$^{-1}$ \\
        $\vect{u}$          & Mobility-free Darcy flux  & m$^2\!\cdot$Pa/s \\
        $u_k$               & Interface mobility-free Darcy flux
        & m$^2\!\cdot$Pa/s \\
        $\vect{w}$          & Advective enthalpy flux      & W/m$^2$ \\
        $w_k$               & Interface advective enthalpy flux
        & W/m$^2$ \\
        $\vect{q}$          & Conductive heat flux         & W/m$^2$ \\
        $q_k$               & Interface conductive heat flux
        & W/m$^2$ \\
        \midrule
        \multicolumn{3}{l}{\textit{Source/sink terms}} \\
        $q^\xi$             & Component mass source/sink   & kg$\cdot$m$^{-3}\!\cdot$s$^{-1}$ \\
        $q^E$               & Energy source/sink           & W/m$^3$ \\
        \midrule
        \multicolumn{3}{l}{\textit{Mixed-dimensional geometry}} \\
        $\Omega_d^i$        & Subdomain $i$ of dim. $d$   & -- \\
        $\Gamma_k$          & Interface between adjacent subdomains   & -- \\
        $\nu_d$             & Specific volume of subdomain with dim. d     & m$^{n-d}$ \\
        $\nu_k$             & Interface specific volume
        & m$^{n-d}$ \\
        $a$                 & Fracture aperture                 & m \\
        $\Xi_k^i$           & Projection: $\Gamma_k \to \Omega_d^i$
        & -- \\
        $\Pi_k^h$           & Trace: $\Omega_h \to \Gamma_k$   & -- \\
        $\hat{\mathcal{S}}_i$ & Set of interfaces to higher-dim.\ subdomains
        & -- \\
        $\check{\mathcal{S}}_i$ & Set of interfaces to lower-dim.\ subdomains
        & -- \\
        $\mathcal{I}_d$     & Index set for subdomains of dim.\ $d$
        & -- \\
        $\tens{I}_n$        & Identity matrix of size $n$       & -- \\
        \bottomrule
    \end{tabular}
\end{table}

\section{Numerical solution}
\label{sec:numerical_solution}
\begin{figure}[H]
    \centering
    \includegraphics[width=0.43\linewidth]{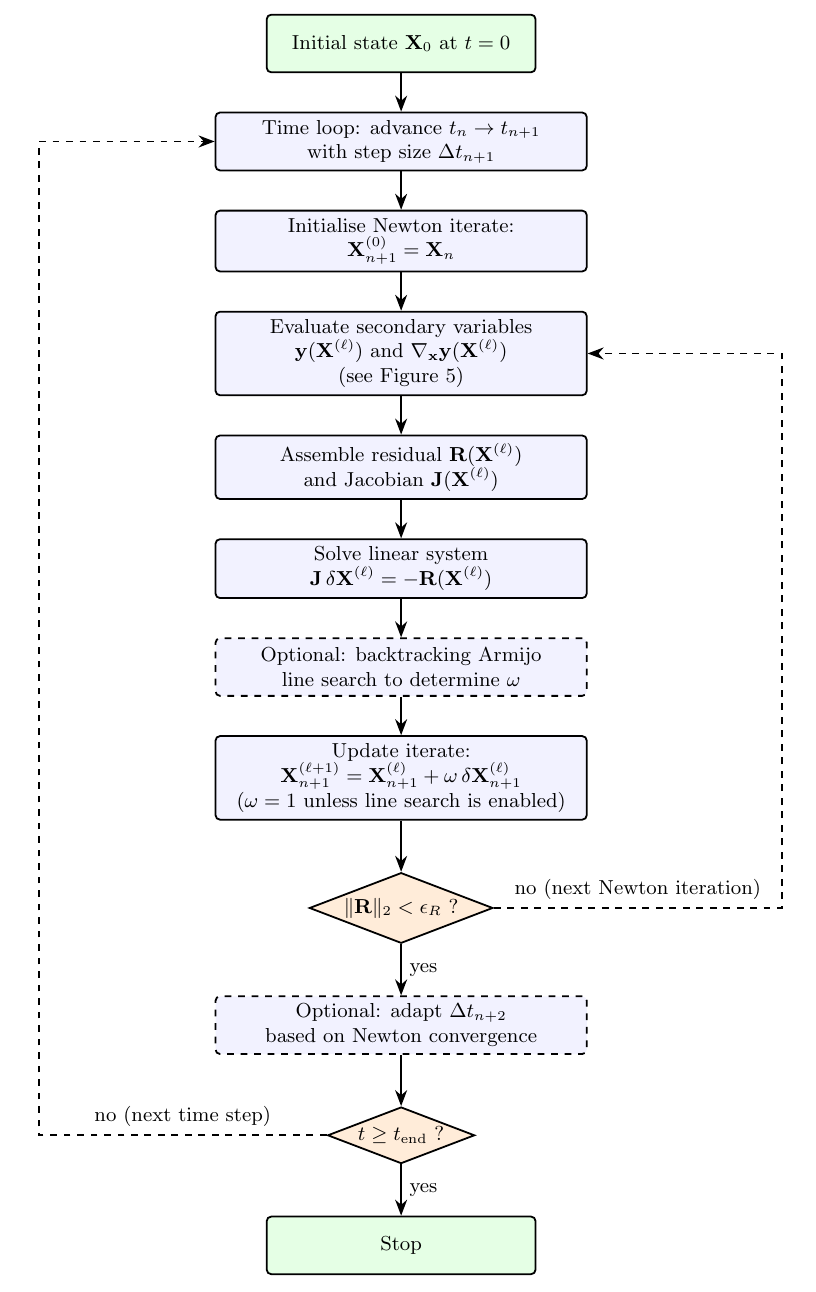}
   \caption{Schematic of the numerical solution loop. Dashed boxes indicate optional components.}
\label{fig:solver_loop}
\end{figure}

This section describes the numerical scheme used to solve the coupled mixed-dimensional system \eqref{eq:component_mass_balance}--\eqref{eq:energy_balance} together with the constitutive relations of Section \ref{constitutive_relations} and the initial and boundary conditions of Section~\ref{sec:ic_bc}. We present the spatial and temporal discretisation, the nonlinear solution strategy, and an efficient evaluation of the secondary map \eqref{eq:secondary_recovery} involving H$_2$O-NaCl fluid properties and state-dependent material properties. The numerical implementation is carried out in the open-source PorePy library \citep{keilegavlen2021porepy}. Figure~\ref{fig:solver_loop} summarises the overall solution procedure.

\subsection{Discretisation scheme}
The three balance equations~\eqref{eq:component_mass_balance},~\eqref{eq:total_mass_balance}, and~\eqref{eq:energy_balance} share a common conservative structure on each subdomain $\Omega_d^i$. We write them in the unified form 
\begin{equation}
	\pd{}{t}\bigl(\nu_d\, \mathcal{A}_i\bigr)
	+ \nabla \cdot \bigl(\nu_d\, \boldsymbol{\mathcal{F}}_i\bigr)
	- \sum_{k \in \hat{\mathcal{S}}_i} 
	\Xi_k^i\, \nu_k\, \mathcal{F}_k
	= \nu_d\, \mathcal{Q}_i,
	\label{eq:generic_balance}
\end{equation}
where $\mathcal{A}_i$ denotes the volumetric accumulation, 
$\boldsymbol{\mathcal{F}}_i$ the subdomain flux, $\mathcal{F}_k$ the 
interface flux, and $\mathcal{Q}_i$ the volumetric source or sink, as given in Table~\ref{tab:generic_balance_adaptation}.
\begin{table}[H]
	\centering
	\caption{Terms of the unified Equation~\eqref{eq:generic_balance} for each balance equation in Section~\ref{sec:model_equations}.}
	\label{tab:generic_balance_adaptation}
	\begin{tabular}{lcccc}
		\toprule
		Balance equation & $\mathcal{A}_i$ & $\boldsymbol{\mathcal{F}}_i$ 
		& $\mathcal{F}_k$ & $\mathcal{Q}_i$ \\
		\midrule
		Component mass~\eqref{eq:component_mass_balance} 
		& $\phi_i \rho_i z_i^\xi$ 
		& $f_i^\xi \lambda_i\, \vect{u}_i$
		& $f_k^\xi \lambda_k\, u_k$
		& $q_i^\xi$ \\[2pt]
		Total mass~\eqref{eq:total_mass_balance} 
		& $\phi_i \rho_i$
		& $\lambda_i\, \vect{u}_i$
		& $\lambda_k\, u_k$
		& $q_i$ \\[2pt]
		Energy~\eqref{eq:energy_balance}
		& $E_i$
		& $\vect{w}_i + \vect{q}_i$
		& $w_k + q_k$
		& $q_i^E$ \\
		\bottomrule
	\end{tabular}
\end{table}

\subsubsection{Spatial discretisation}
\label{sec:spatial_discretization}
A cell-centered finite volume scheme is used. Each subdomain $\Omega_d^i$ is partitioned into non-overlapping cells $\mathcal{T}_d^i = \{\tau\}$, where $\tau$ denotes a generic cell with volume $|\tau|$ and faces $\sigma \subset \partial \tau$. The mesh is conforming across matrix-fracture and fracture-intersection interfaces \citep{keilegavlen2021porepy}. The primary variables $\vect{x}_\tau = (p_\tau, h_\tau, z_\tau)$ are stored at cell centers. We denote by $\mathbf{n}_{\tau,\sigma}$ the outward unit normal to $\sigma$ on $\tau$, and by $\nu_{d,\tau}$ the specific volume evaluated cell-wise on $\tau$. Integrating~\eqref{eq:generic_balance} over $\tau$ and applying the divergence theorem yields the semi-discrete form:
\begin{equation}
    |\tau|\, \frac{d}{dt}\bigl(\nu_{d,\tau}\, \mathcal{A}_\tau\bigr)
    + \nu_{d,\tau}\sum_{\sigma \subset \partial \tau} \mathcal{F}_{\tau,\sigma}
    - \sum_{k \in \hat{\mathcal{S}}_i} 
      \bigl(\Xi_k^i\, \nu_k\, \mathcal{F}_k\bigr)_\tau
    = |\tau|\, \nu_{d,\tau}\, \mathcal{Q}_\tau,
    \label{eq:fv_generic}
\end{equation}
where $\mathcal{A}_\tau$ and $\mathcal{Q}_\tau$ are the cell-averaged values of $\mathcal{A}_i$ and $\mathcal{Q}_i$ on $\tau$. The flux through face $\sigma$ is given by
\begin{equation}
    \mathcal{F}_{\tau,\sigma} 
    := \int_\sigma \boldsymbol{\mathcal{F}}_i \cdot 
    \mathbf{n}_{\tau,\sigma}\,\mathrm{d}s.
    \label{eq:face_flux}
\end{equation}
The subdomain flux $\boldsymbol{\mathcal{F}}_i$ (see Table~\ref{tab:generic_balance_adaptation}) consists of an advective contribution $\mathcal{W}_i\, \vect{u}_i,$ where  $\mathcal{W}_i\in \{\lambda_i,\, f_i^\xi \lambda_i,\, 
\sum_\gamma \lambda_i^\gamma h_i^\gamma\}$ is the cell-wise nonlinear advective weight, with the choice determined by the balance equation, and for the energy balance only, an additional conductive contribution $\vect{q}_i.$ Both $\vect{u}_i$ and $\vect{q}_i$ are discretised by multi-point flux approximation (MPFA)~\citep{aavatsmark2002introduction, schneider2020coupling}: for an interior face $\sigma$,
\begin{equation}
    U_\sigma 
    := \int_\sigma \vect{u}_i \cdot \mathbf{n}_{\tau,\sigma}\,\mathrm{d}s
    = -\sum_{M \in \mathcal{S}_\sigma} t_{\sigma,M}^p\, p_M,
    \qquad
    Q_\sigma 
    := \int_\sigma \vect{q}_i \cdot \mathbf{n}_{\tau,\sigma}\,\mathrm{d}s
    = -\sum_{M \in \mathcal{S}_\sigma} t_{\sigma,M}^T\, T_M,
    \label{eq:mpfa}
\end{equation}
where $\mathcal{S}_\sigma$ is the MPFA stencil, i.e., the set of cells $M$ contributing to the flux through $\sigma,$ and $\{t_{\sigma,M}^p\}, \{t_{\sigma,M}^T\}$ are the corresponding local transmissibilities computed as in~\citet{aavatsmark2002introduction}; $p_M$ and $T_M$ are the cell pressure and temperature. For boundary faces, the stencil also incorporates the prescribed Dirichlet or Neumann data of Section~\ref{sec:ic_bc}. 

The nonlinear weight $\mathcal{W}_i$ is upwinded with respect to the 
sign of $U_\sigma$:
\begin{equation}
    \widehat{\mathcal{W}}_{i,\sigma} = 
    \begin{cases}
        \mathcal{W}_i\bigl(\vect{x}_\tau\bigr), & U_\sigma \geq 0,\\[2pt]
        \mathcal{W}_i\bigl(\vect{x}_{\tau'}\bigr), & U_\sigma < 0,
    \end{cases}
    \label{eq:upwind}
\end{equation}
where $\tau'$ denotes the cell on the opposite side of $\sigma$ from $\tau$. The discrete face flux is then
\begin{equation}
    \mathcal{F}_{\tau,\sigma} 
    = \widehat{\mathcal{W}}_{i,\sigma}\, U_\sigma
    \label{eq:face_mass_flux}
\end{equation}
for the mass fluxes, and
\begin{equation}
 \mathcal{F}_{\tau,\sigma} 
= \widehat{\mathcal{W}}_{i,\sigma}\, U_\sigma + Q_\sigma
\label{eq:face_energy_flux}
\end{equation}
for the energy flux. The projected interface flux contribution in~\eqref{eq:fv_generic} is defined as 
\begin{equation}
\bigl(\Xi_k^i\, \nu_k\, \mathcal{F}_k\bigr)_\tau 
:= \int_\tau \Xi_k^i\, \nu_k\, \mathcal{F}_k\,\mathrm{d}V
= \begin{cases}
    |\sigma_k|\, \nu_k\, \mathcal{F}_k & \text{if } \sigma_k \subset \partial\tau,\\
    0 & \text{otherwise},
\end{cases}
\end{equation}
where $|\sigma_k|$ is the area of the face shared between $\tau$ and $\Gamma_k$. The interface flux $\mathcal{F}_k$ (see Table~\ref{tab:generic_balance_adaptation}) is obtained analogously to~\eqref{eq:face_mass_flux} and~\eqref{eq:face_energy_flux}: 
$\mathcal{F}_k = \widehat{\mathcal{W}}_k\, U_k$ for the interface mass fluxes and $\mathcal{F}_k = \widehat{\mathcal{W}}_k\, U_k + Q_k$ for the interface energy flux. Here, $U_k$ and $Q_k$ are obtained from~\eqref{eq:interface_darcy} and~\eqref{eq:interface_energy_fluxes}, respectively, and the interface upwind weight $\widehat{\mathcal{W}}_k$ is taken from $\tau$ (the higher-dimensional cell) if $U_k \geq 0$ and from the lower-dimensional neighbour otherwise.

\subsubsection{Temporal discretisation}
\label{subsec:temporal_discretization}
We seek the solution of~\eqref{eq:fv_generic} over a simulation interval $[0, t_{\mathrm{end}}]$. Collecting the cell-centered primary variables $\vect{x}_\tau =(p_\tau, h_\tau, z_\tau)$ over all cells of all subdomains into a global state vector 
$\vect{X}(t)$, the semi-discrete system~\eqref{eq:fv_generic} takes the form of a nonlinear system of ordinary differential equations in time
\begin{equation}
    \frac{\mathrm{d}}{\mathrm{d}t} 
    \mathbf{M}\bigl(\mathbf{X}(t)\bigr)
    + \mathbf{G}\bigl(\mathbf{X}(t)\bigr) = \mathbf{0},
    \label{eq:semi_discrete_ode}
\end{equation}
where $\vect{M}(\vect{X})$ collects the cell-wise accumulation contributions and $\mathbf{G}(\vect{X})$ collects all flux and source terms in~\eqref{eq:fv_generic}; both are nonlinear functions of $\vect{X}$ through the constitutive relations of Section~\ref{constitutive_relations}.

We integrate~\eqref{eq:semi_discrete_ode} by the fully implicit backward Euler scheme on a sequence of time levels $0 = t_0 < t_1 < \cdots < t_N = t_{\mathrm{end}}$ with step sizes $\Delta t_{n+1} = t_{n+1} - t_n,$ which may be constant or adapted. The resulting discrete nonlinear system reads
\begin{equation}
    \frac{\mathbf{M}\bigl(\mathbf{X}_{n+1}\bigr) 
        - \mathbf{M}\bigl(\mathbf{X}_{n}\bigr)}{\Delta t_{n+1}}
    + \mathbf{G}\bigl(\mathbf{X}_{n+1}\bigr) = \mathbf{0},
    \label{eq:backward_euler}
\end{equation}
where $\mathbf{X}_n \approx \mathbf{X}(t_n)$, $\mathbf{M}(\mathbf{X}_n)$ is taken from the previous time step $t_n$, and all remaining terms are evaluated at $t_{n+1}$. When adaptive time stepping is used, $\Delta t_{n+1}$ is adjusted based on the convergence behaviour of the nonlinear solver for \eqref{eq:backward_euler}, as described in Section~\ref{sec:nonlinear_linear_solver}. This corresponds to the outer ``Time loop'' in Figure~\ref{fig:solver_loop}.

\subsection{Nonlinear solution strategy}
\label{sec:nonlinear_linear_solver}
Each step of the time loop in Figure~\ref{fig:solver_loop} invokes Newton's method on the nonlinear discrete system~\eqref{eq:backward_euler}. The nonlinear residual associated with \eqref{eq:backward_euler} is defined as
\begin{equation}
	\mathbf{R}(\mathbf{X}_{n+1}) :=
	\frac{\mathbf{M}(\mathbf{X}_{n+1}) - \mathbf{M}(\mathbf{X}_n)}{\Delta t_{\,n+1}}
	+ \mathbf{G}(\mathbf{X}_{n+1}).
	\label{eq:residual}
\end{equation}
At Newton iteration $\ell$, the linear system reads
\begin{equation}
	\mathbf{J}\bigl(\mathbf{X}_{n+1}^{(\ell)}\bigr)
	\, \delta \mathbf{X}_{n+1}^{(\ell)}
	=
	-\mathbf{R}\bigl(\mathbf{X}_{n+1}^{(\ell)}\bigr),
	\label{eq:newton}
\end{equation}
where
\begin{equation}
	\delta \mathbf{X}_{n+1}^{(\ell)}
	:=
	\mathbf{X}_{n+1}^{(\ell+1)} - \mathbf{X}_{n+1}^{(\ell)}
    \label{eq:newton_update}
\end{equation}
is the Newton update and $\mathbf{J} = \partial \mathbf{R}/\partial \mathbf{X}$ is the Jacobian assembled by automatic differentiation in PorePy.

The persistent-variable formulation~\eqref{var:primary_variables} preserves $\mathbf{X}$ across phase transitions, eliminating the need for variable switching. The linear system~\eqref{eq:newton} is solved for $\delta \mathbf{X}_{n+1}^{(\ell)}$ using a sparse direct solver, which is sufficient for the 1D and 2D problems considered in this work.

The Convergence of the Newton iteration is declared when
\begin{equation}
	\|\mathbf{R}(\mathbf{X}_{n+1}^{(\ell)})\|_2 < \epsilon_R,
	\label{eq:convergence}
\end{equation}
where $\|\cdot\|_2$ denotes a $l_2$-norm and $\epsilon_R$ is a prescribed tolerance; 
specific values are given in 
Sections~\ref{sec:verification} and~\ref{sec:application}. 

The time step in \eqref{eq:backward_euler} is adapted as follows: if convergence is achieved within a target number of iterations $\ell_{\mathrm{target}}$, the next step size is increased by a factor $\alpha > 1$; if convergence requires more than $\ell_{\mathrm{target}}$ iterations, the step size is kept unchanged; if convergence fails within a maximum iteration count $\ell_{\max}$, the step is rejected, $\Delta t$ is reduced by a factor $\beta < 1$, and the step is re-attempted.

\subsubsection{Backtracking Armijo line search}
\label{sec:armijo_line_search}
For problems with strong nonlinearities, particularly near phase boundaries, where the thermodynamic closure \eqref{eq:secondary_recovery} is non-smooth, the standard Newton update from Equation \eqref{eq:newton_update}:
\begin{equation}
	\mathbf{X}_{n+1}^{(\ell + 1)} = \mathbf{X}_{n+1}^{(\ell)} + \delta\mathbf{X}_{n+1}^{(\ell)}
\end{equation}
can overshoot and stall convergence or push the state into unphysical regions. To improve robustness, we optionally enable a backtracking Armijo line search:
\begin{equation}
    \vect{X}^{(\ell+1)}_{n+1} = \vect{X}^{(\ell)}_{n+1}
    + \omega\,\delta\vect{X}^{(\ell)}_{n+1},
\end{equation}
where $\omega\in (0,1]$ is the step size determined by the backtracking condition \citep{nocedal2006numerical}. Starting from $\omega = 1$, the step size is reduced by a factor
$\overline{c} \in (0,1)$ until the sufficient decrease condition is satisfied:
\begin{equation}
    \Phi(\vect{X}^{(\ell)}_{n+1}
    + \omega\,\delta\vect{X}^{(\ell)}_{n+1})
    \leq (1 - 2\theta\omega)\,
    \Phi(\vect{X}^{(\ell)}_{n+1}).
\end{equation}
Here, $\Phi(\vect{X}) = \tfrac{1}{2}\|\vect{R}(\vect{X})\|_2^2$ is the residual objective function with $\vect{R}$ defined by \eqref{eq:residual}, and $\theta > 0$ is the Armijo parameter. When line search is disabled, $\omega = 1$ is used. 

\subsubsection{Evaluation of secondary variables}
\label{subsec:fluid_properties_obl}
\begin{figure}[H]
    \centering
    \includegraphics[width=0.5\linewidth]{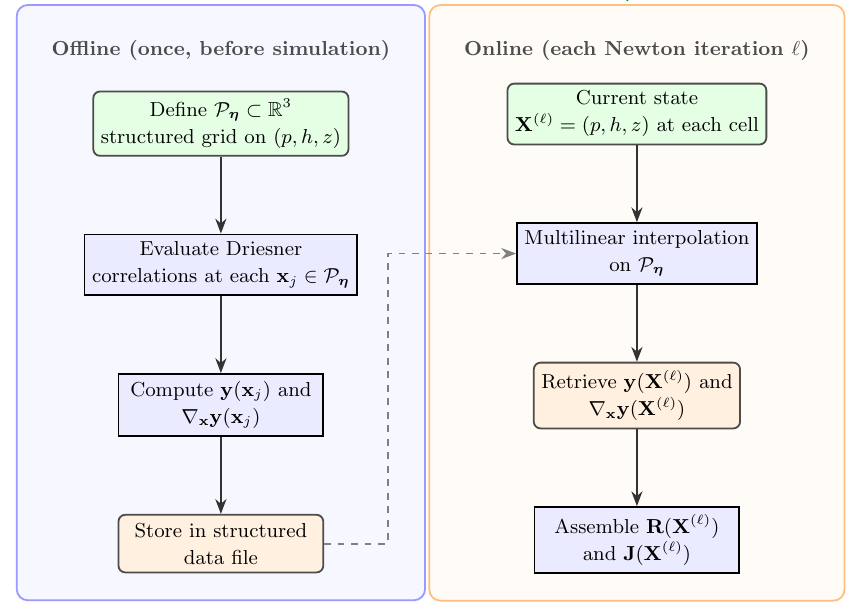}
   \caption{Workflow for evaluating the secondary-variable map \eqref{eq:secondary_recovery}. Offline (left): the secondary-variable map $\mathbf{y}$ and its gradient $\nabla_{\mathbf{x}} \mathbf{y}$ are precomputed from the Driesner correlations on a structured grid $\mathcal{P}_{\boldsymbol{\eta}}$ over the primary-variable space $(p, h, z)$ and stored in a data file. Online (right): at each Newton iteration, $\mathbf{y}$ and $\nabla_{\mathbf{x}} \mathbf{y}$ are retrieved at the current state by multilinear interpolation on $\mathcal{P}_{\boldsymbol{\eta}}$ and used to assemble the residual $\mathbf{R}$ and Jacobian $\mathbf{J}$.}
   \label{fig:obl_workflow}
\end{figure}

The assembly of the residual $\vect{R}$ and the Jacobian $\vect{J}$ in Equation \eqref{eq:newton} requires the nonlinear secondary variable map \eqref{eq:secondary_recovery} and its gradient for every grid cell at each Newton iteration. This step usually dominates the total computational cost. To reduce this cost, we adopt an approach inspired by the operator-based linearisation (OBL) technique of \citet{voskov2017operator},  illustrated schematically in Figure~\ref{fig:obl_workflow}.

Let $\mathcal{P}_{\boldsymbol{\eta}} \subset \mathbb{R}^3$ be a structured grid on the primary variable space $\mathbf{x} :=(p, h, z)$ with mesh size $\boldsymbol{\eta}=(\eta_p, \eta_h, \eta_z)$. The
map $\mathbf{y}(\mathbf{x}_j)$ and its gradient $\nabla_{\mathbf{x}} \mathbf{y}(\mathbf{x}_j)$ are precomputed at every node $\mathbf{x}_j \in \mathcal{P}_{\boldsymbol{\eta}}$ using the saltwater correlations of~\citet{driesner2007system}. For any arbitrary state $\mathbf{x}$ the secondary variables are recovered via multilinear interpolation on $\mathcal{P}_{\boldsymbol{\eta}}$. The linearisation error is controlled directly by ${\boldsymbol{\eta}}$ \citep{voskov2017operator}.

The choice of ${\boldsymbol{\eta}}$ involves a trade-off between interpolation accuracy and computational cost. A finer grid reduces the linearisation error but increases both the storage memory for precomputed values and the offline computation time. In this study, ${\boldsymbol{\eta}}$ is selected heuristically to resolve the strong nonlinearities of the H$_2$O-NaCl phase diagram, particularly the sharp gradients in the secondary variables near phase transition boundaries. The chosen resolution, with mesh sizes $\eta_p = \SI{10}{\bar}$, $\eta_h = \SI{10}{\kilo\joule\per\kilo\gram}$, and $\eta_z = 0.01$, is found to reproduce the benchmark of Section \ref{sec:verification} in close agreement with the CSMP++ reference solution.   

The quantities evaluated through this approximation include temperature, phase saturations, compositions, densities, viscosities, and enthalpies. Dependent algebraic relations such as \eqref{eq:porosity_update}-\eqref{eq:permeability_kozeny_carman} are computed accordingly.

\section{Numerical benchmarking}
\label{sec:verification}
We verify the numerical framework of Section \ref{sec:numerical_solution} implemented in PorePy against the one-dimensional salt-dissolution benchmark proposed by \citet{weis2014hydrothermal} for code comparison with the CSMP++ simulator. Despite its simple geometry, the benchmark 
tests the main components of the formulation, including phase transitions between liquid, vapour and solid halite within the persistent-variable formulation; the nonlinear coupling between primary variables; and the feedback of halite dissolution on the hydraulic properties of the porous medium. Although both simulators utilise a persistent-variable approach to avoid explicit variable switching, they differ in their numerical implementation. CSMP++ uses a control volume finite element method with a semi-implicit sequential time-stepping scheme \citep{weis2014hydrothermal}. In contrast, our framework adopts a cell-centered finite volume method with MPFA and upwind flux discretisation and a fully implicit time-stepping scheme (see Section \ref{sec:numerical_solution}). While this one-dimensional benchmark does not test the mixed-dimensional coupling, it provides a necessary foundation for the fractured-domain applications in Section \ref{sec:application}. The simulation scripts to reproduce the results in Sections \ref{sec:verification}-\ref{sec:application} are available in a Docker container \citep{oguntola2026docker}.

\paragraph{Problem setup.}
We consider the H$_2$O--NaCl system as a two-component mixture comprising three possible phases: liquid, vapour, and solid halite. The simulation domain is a horizontal one-dimensional column $\Omega = (0, L)$ with $L = 2$~km, discretised with a uniform mesh size $\Delta \zeta = 10$~m. Gravity and capillary effects are neglected, and the lateral boundaries of the (pseudo-one-dimensional) domain are impermeable and adiabatic. 

Initially, the domain contains halite-saturated brine in local equilibrium with precipitated halite at a uniform initial saturation $s_0^\mathrm{hal} = 0.1$ and temperature $T_0 = 150^\circ$C. The initial pressure $p_0$ is distributed linearly between the boundary values. Given that our formulation employs persistent primary variables \eqref{var:primary_variables}, the initial overall salt mass fraction is determined pointwise by solving
\begin{equation}
z_{0}(\zeta) = \arg\min_z \bigl|\tilde{s}^{\text{hal}}\bigl(p^0(\zeta), T^0, z\bigr) - s^{\mathrm{hal}}_0\bigr|,
\label{eq:dissol-initial-salt-fraction}
\end{equation}
where $\tilde{s}^{\text{hal}}$ is the halite saturation evaluated from the saltwater EOS \citep{driesner2007system}. The initial specific enthalpy then follows as $h_0(\zeta) = \tilde{h}\bigl(p_0(\zeta), T_0, z_{0}(\zeta)\bigr)$.

Dirichlet boundary conditions are imposed at both ends of the domain. At the inlet ($\zeta = 0$), single-phase water vapour is prescribed with $p_{\mathrm{in}} = 4$~MPa, $T_{\mathrm{in}} = 300^\circ$C, and $z_{\mathrm{in}} = 0$. At the outlet ($\zeta = L$), a halite-saturated state consistent with the initial condition is maintained, with $p_{\mathrm{out}} = 1$~MPa, $T_{\mathrm{out}} = 150^\circ$C. The boundary enthalpies $h_{\mathrm{in}}$ and $h_{\mathrm{out}}$ are recovered from the closure map \eqref{eq:secondary_recovery}, and the outlet salt $z_{\mathrm{out}}$ is obtained from \eqref{eq:dissol-initial-salt-fraction}. The resulting pressure gradient drives hot vapour from the inlet into the halite-saturated brine, progressively heating the domain and dissolving the initially contained halite.

\paragraph{Physical and numerical parameters.}
The rock matrix properties are homogeneous with a reference permeability $K^0 = 10^{-15}\,\mathbf{I}$~m$^2$, porosity $\phi_0 = 0.1$, rock density $\rho_r = 2700$~kg/m$^3$, specific heat capacity $c_{p,r} = 880$~J/(kg$\cdot$K), and thermal conductivity $\kappa_r = 2$~W/(m$\cdot$K), consistent with values from \cite{weis2014hydrothermal}. The relative permeabilities for the mobile phases $\gamma$ follow a Brooks-Corey form with residuals \citep{oguntola2020robust, falko2021brine},

\begin{equation}
k_{r\gamma}(s^\gamma) = \begin{cases} 0, & s^\gamma \le R_\gamma, \\[2pt] \dfrac{s^\gamma - R_\gamma}{1 - (R_{\mathrm{liq}} + R_{\mathrm{vap}})}, & s^\gamma > R_\gamma, \end{cases}
\label{eq:multiphase_relperm}
\end{equation}
where $R_{\mathrm{liq}} = 0.3$ and $R_{\mathrm{vap}} = 0$ are the residual liquid and vapour saturations. The simulation is advanced to $t_{\mathrm{end}} = 2000$~years using a fixed time-step size $\Delta t = 365$~days, tolerance $\epsilon_R = 1.0e-5,$ and disabled line search method. 

\begin{figure}[htp!]
    \centering
    \begin{minipage}{0.48\linewidth}
        \centering
        \includegraphics[width=\linewidth]{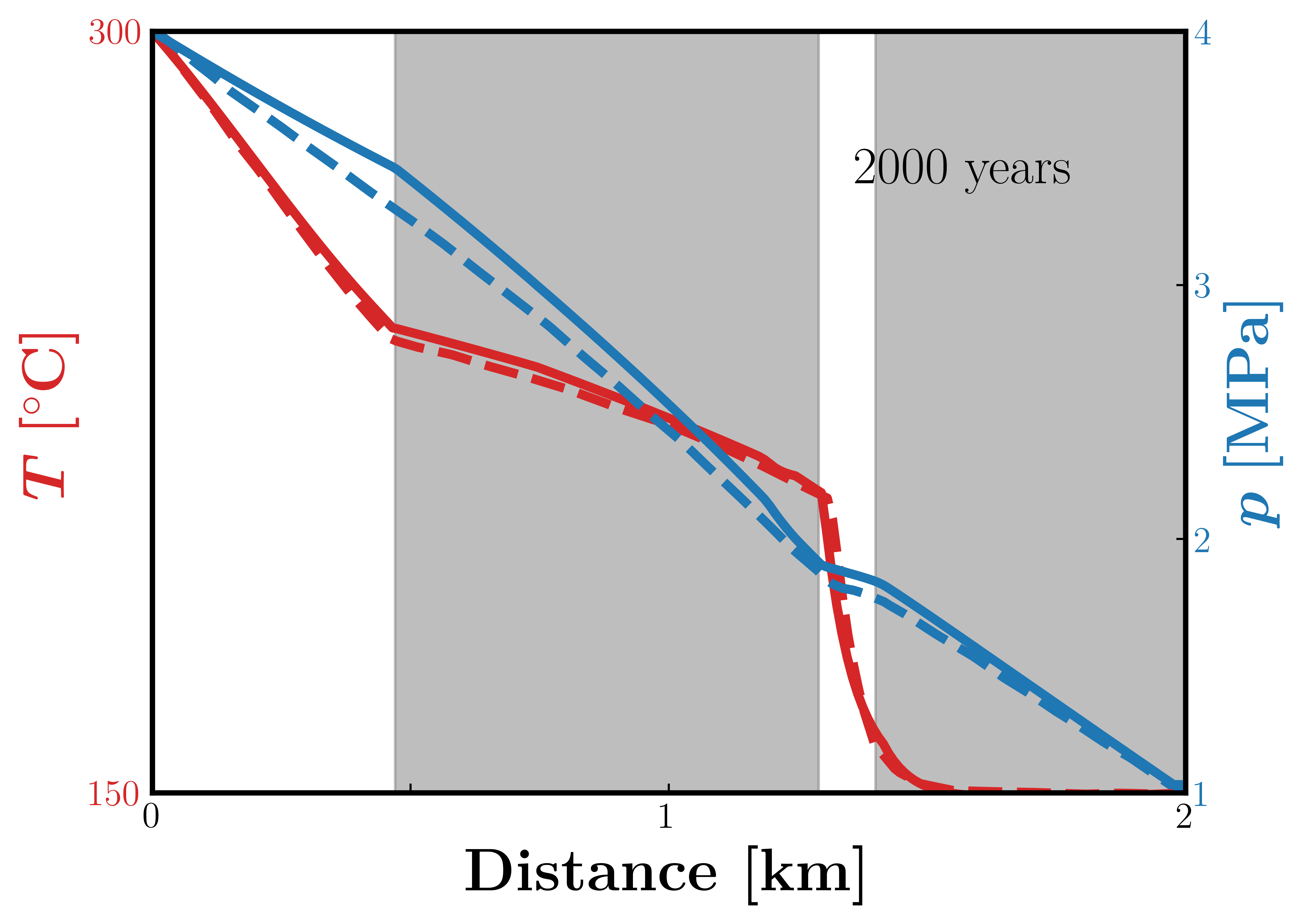}
    \end{minipage}
    \hfill
    \begin{minipage}{0.48\linewidth}
        \centering
        \includegraphics[width=\linewidth]{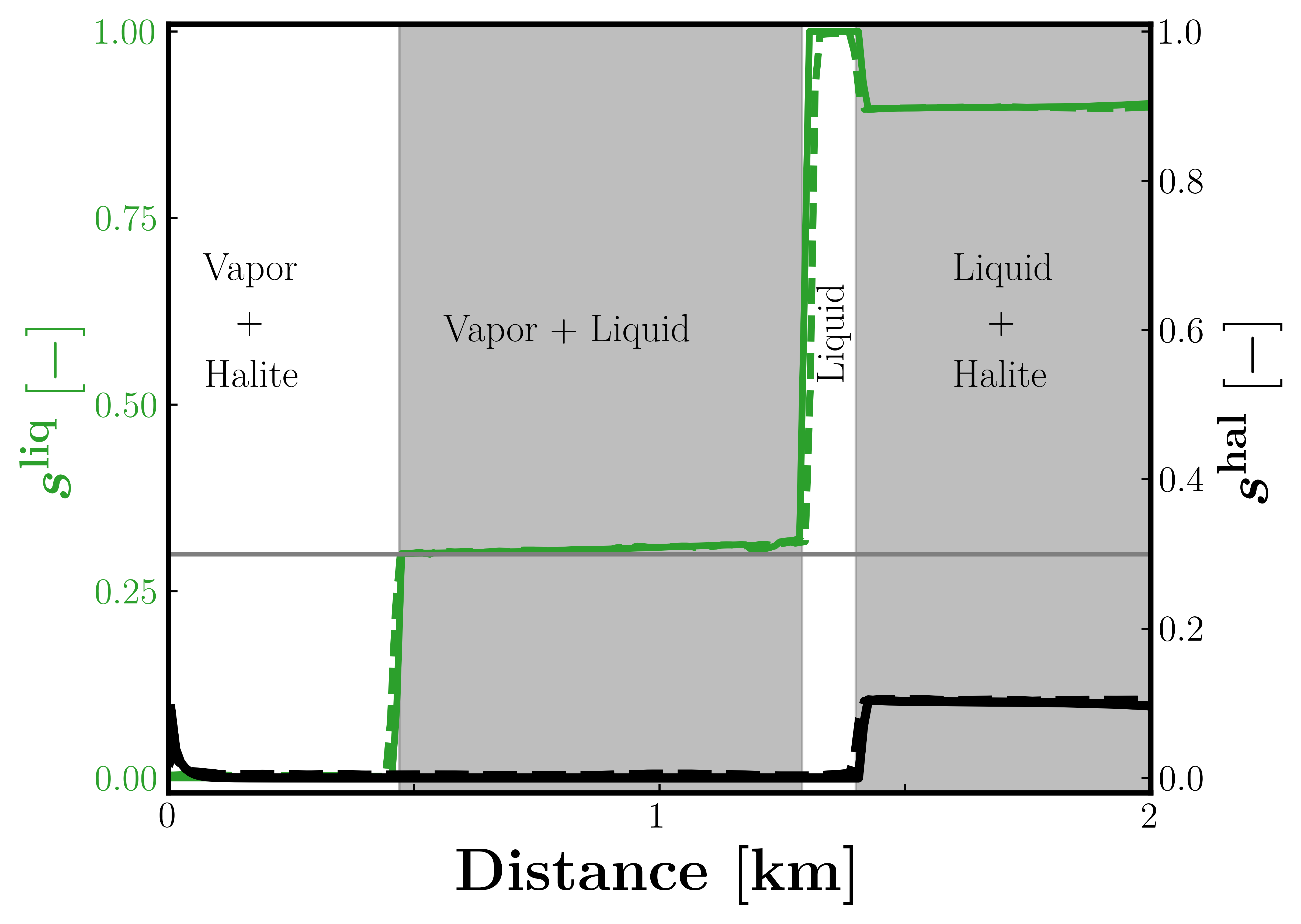}
    \end{minipage}
    \caption{One-dimensional salt-dissolution benchmark at time $t = 2000$~years. Left: fluid pressure (blue) and temperature (red). Right: liquid saturation (green), halite saturation (black), and residual liquid saturation (horizontal gray line). Solid curves are the PorePy results; dashed curves are the CSMP++ reference of~\citet{weis2014hydrothermal}, digitized from the published figure.}
    \label{fig:benchmark_solution}
\end{figure}

\paragraph{Results.} Figure \ref{fig:benchmark_solution} compares the profiles of pressure, temperature, and liquid and halite saturations at time $t=2000$~years from the PorePy implementation (solid curves) with the CSMP++ reference solution of \citet{weis2014hydrothermal} (dashed curves, digitised from their Figures 6C and 6D). Four distinct phase regions are observed along the flow direction, and their order, extent, and transition locations agree closely between the two solutions. 

Near the inlet, a vapour+halite region extends to roughly $\zeta \approx 0.48$~km, where the high-enthalpy hot vapour has displaced the mobile brine and dissolved most of the solid halite initially present; in this region, liquid saturation is essentially zero and only a small residual halite saturation remains. Downstream, a vapour+liquid boiling zone occupies approximately $0.48 \lesssim \zeta \lesssim 1.6$~km, in which the liquid saturation settles at the residual value $s^{\mathrm{liq}} = R_{\mathrm{liq}} = 0.3$ and the temperature follows the H$_2$O--NaCl two-phase coexistence curve at the local pressure and salinity. A narrow single-phase liquid region near $\zeta \approx 1.6$~km forms, where the thermal energy carried by the advancing fluid is insufficient to induce boiling at the local pressure; here the liquid saturation rises sharply to unity. Beyond this front, a liquid+halite region extends to the outlet and remains near the initial state.

The two solutions agree on the locations of all three transition fronts, on the residual liquid saturation within the boiling zone, and on the temperature and pressure profiles through all four regimes. The most visible difference is in the pressure profile within the boiling zone, where the PorePy solution predicts slightly higher pressures than CSMP++; a smaller discrepancy is also apparent in the shape of the narrow liquid region. These differences are attributed to digitisation artefacts and differences in spatial discretisation (cell-centered vs. node-centered) and fluid-property evaluation between the two codes.

\section{Application}
\label{sec:application}
In this section, we apply the numerical framework to a two-dimensional high-enthalpy fractured geothermal reservoir. This case study demonstrates the model's capability to predict halite precipitation and dissolution patterns and their impact on reservoir performance under production-driven flow in a mixed-dimensional setting. We present three examples. The first two (Sections~\ref{sec:example1} and~\ref{sec:example2}) consider a
reservoir with disconnected fractures and examine the effect of the clogging exponent~$\varphi$ in Equation~\eqref{eq:aperture_update} on near-wellbore halite precipitation and the corresponding impact on energy recovery. The third (Section~\ref{sec:example3}) considers a connected fracture network to investigate the role of fracture connectivity.
\begin{figure}
    \centering
    \begin{minipage}{0.8\linewidth}
        \centering
        \includegraphics[width=\linewidth]{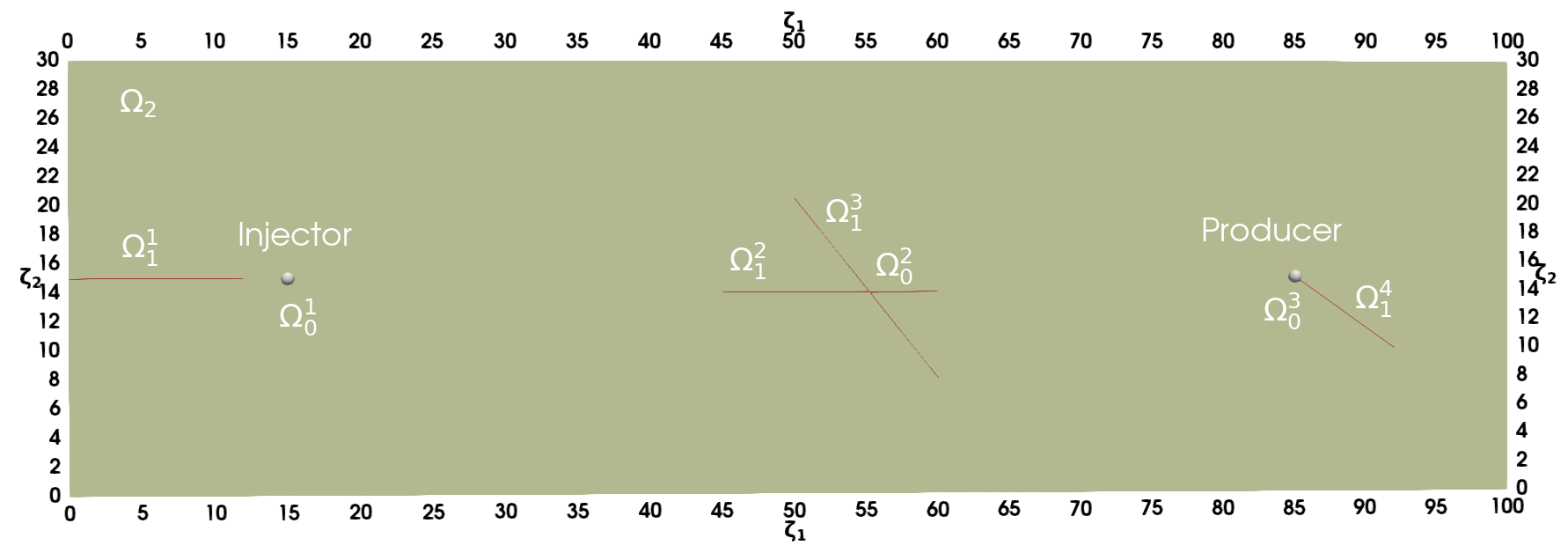}
    \end{minipage}
    \hfill
    \begin{minipage}{0.8\linewidth}
        \centering
        \includegraphics[width=\linewidth]{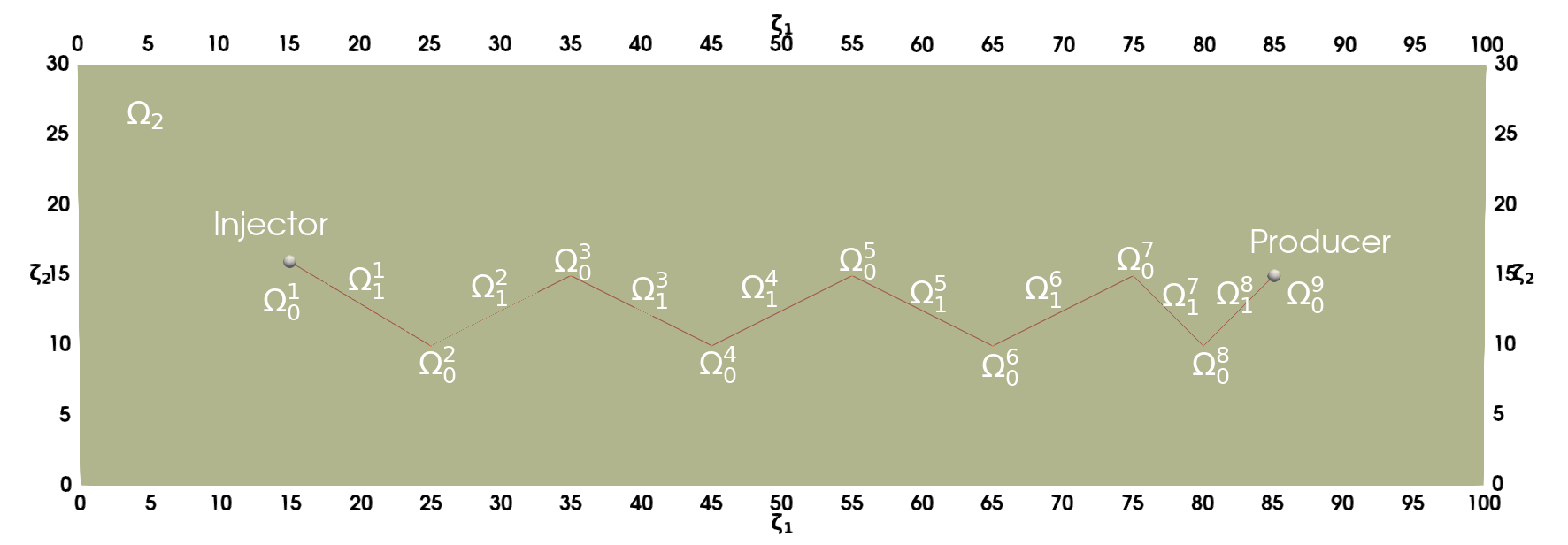}
    \end{minipage}
    \caption{Mixed-dimensional representation of the 2D fractured geothermal reservoir $(100 \times 30 \text{ m})$. The diagrams depict the spatial decomposition of the reservoir: the 2D matrix ($\Omega_2$), 1D fracture segments ($\Omega_1^i$), and 0D subdomains ($\Omega_0^i$) representing fracture intersections and point-wells. The top panel illustrates a network of four disconnected fractures with a single intersection. The bottom panel illustrates a connected chain of eight fractures and seven intersections between the well subdomains. In both cases, the injection well is located at approximately $(15, 15)$~m and the production well is at $(85, 15)$~m. All external boundaries are impermeable and adiabatic, with the fluid flow driven entirely by the pressure differential between the two 0D wells.}
    \label{fig:application_geometry}
\end{figure}
\begin{table}[!ht]
\centering
\caption{Physical and numerical parameters for Examples $1 - 3$}
\label{tab:common_parameters}
\small
\renewcommand{\arraystretch}{1.2}
\begin{tabular}{@{}l l l@{}}
\toprule
\textbf{Parameter} & \textbf{Value} & \textbf{Unit} \\
Reference porosity ($\phi^0$)       & 0.1                         & -- \\
Reference permeability ($\tens{K}^0$) & $10^{-15}\,\tens{I}$         & m$^2$ \\
Rock density ($\rho_r$)             & 2700                        & kg\,m$^{-3}$ \\
Specific heat capacity ($c_{p,r}$)  & 880                         & J\,kg$^{-1}$\,K$^{-1}$ \\
Thermal conductivity ($\kappa_r$)   & 2.0                         & W\,m$^{-1}$\,K$^{-1}$ \\
Normal permeability ($K_k$)         & $10^{-13}$                  & m$^2$ \\
Well radius ($r_w$)                 & 0.1                         & m \\
Well-cell thickness ($h$)                    & 1.0                & m \\
Skin factor ($s$)                            & 0.0                & -- \\ 
Production BHP ($p_{\mathrm{BHP}}$) & 7.0                         & MPa \\
Injection temperature ($T_{\mathrm{inj}}$) & 300.65 ($\approx 27.5\,^\circ$C) & K \\
Injection salt fraction ($z_{\mathrm{inj}}$) & $1.0 \times 10^{-4}$  & -- \\
Residual liquid saturation ($R_{\mathrm{liq}}$) & 0.3& -- \\
Residual vapour saturation ($R_{\mathrm{vap}}$)  & 0.0              & -- \\
Reference aperture ($a^0$) & $10^{-3}$         & m \\
Minimum aperture ($a_{\min}$)       & $10^{-4}$                   & m \\
Backtracking factor ($\overline{c}$)  & $0.8$                   & -- \\
Armijo parameter ($\theta$)  & $10^{-2}$  & -- \\
\bottomrule
\end{tabular}
\end{table}
\subsection{Example 1}
\label{sec:example1}
We consider a rectangular domain $\Omega = [0, 100] \times [0, 30]$\,m representing a horizontal cross-section of a high-enthalpy geothermal reservoir at depth (Figure~\ref{fig:application_geometry}, top panel). The porous matrix $\Omega_2$ contains four embedded fractures $\Omega^1_1, \ldots, \Omega^4_1$, of which $\Omega^2_1$ and $\Omega^3_1$ intersect to form a single 0D intersection subdomain $\Omega^1_0$. The remaining fractures are isolated, so the fracture network does not provide a continuous pathway between the wells; inter-well flow is therefore carried primarily by the matrix. Injection and production wells, represented as 0D point grids $\Omega^1_0$ and $\Omega^3_0$, are placed at $(15,\,15)$\,m and $(85,\,15)$\,m, respectively, and coupled to the matrix subdomain $\Omega_2$ via well--matrix interfaces. The domain is discretised with an unstructured simplex mesh with matrix cell size 0.7\,m and fracture cell size 1.0\,m.

The reservoir fluid is modelled as a binary H$_2$O--NaCl system that can coexist in up to three phases: liquid brine, vapour, and solid halite. Fluid and rock properties are given in Table \ref{tab:common_parameters}. Initially, the reservoir is in a two-phase state consisting of liquid brine and solid halite under uniform conditions: pressure $p_0 = 10.5$\,MPa, specific enthalpy $h_0 = 1002$\,kJ/kg, overall salt mass fraction $z_0 = 0.40.$ This thermodynamic state yields an initial temperature of $T_0 = 586.65$\,K ($\approx 313.5\,^\circ$C) and an initial halite saturation of $s_0^{\text{hal}} = 0.01055$. The fracture aperture clogging exponent in Equation \eqref{eq:aperture_update} is set to $\varphi = 0.1,$ to allow for mild clogging by precipitated halite in fractures, especially near the production well. The phase relative permeabilities are given by Equation~\eqref{eq:multiphase_relperm}, with the residual liquid saturation scaled by the mobile pore-space fraction
($1 - s^{\text{hal}})$ to account for pore-volume reduction by precipitated halite. All external boundaries are impermeable and adiabatic ($Q^\xi_N = 0$, $q^E_N = 0$; cf. Equation~\eqref{eq:bc_neu}), so that flow is driven entirely by the well pair ($\Omega_0^1, \Omega_0^3$).

The injection well $\Omega_0^1$ injects cold, low-salinity water, defined by a constant volumetric mass rate ${q}_{\text{inj}} = 0.28$~kgm$^{-3}$s$^{-1}$. The quantity ${q}_{\text{inj}}$ enters the total mass balance \eqref{eq:total_mass_balance} as a source term, and the corresponding component source for each $\xi$ is given by $q^{\xi}_\text{inj}=z^\xi_\text{inj} {q}_{\text{inj}},$ where $z^\xi_\text{inj}$ is the injected overall mass fraction of component $\xi.$ These component sources enter the component mass balance~\eqref{eq:component_mass_balance}. The standard energy balance \eqref{eq:energy_balance} is replaced by a Dirichlet temperature constraint,
\begin{equation}
    T = T_{\text{inj}}\qquad\text{on }\Omega_0^1,
\end{equation}
where $T_{\text{inj}}$ is the temperature of the injected fluid (see Table \ref{tab:common_parameters} for fluid properties). The injected fluid specific enthalpy is enforced by the temperature-primary variable relation of $\eqref{eq:secondary_recovery}.$

The production well $\Omega^3_0$ operates at a fixed bottom-hole pressure $p_{\text{BHP}}.$ The resulting production rate $q_{\text{prod}},$ given by
\begin{equation}
q_{\text{prod}} = -\lambda\,\text{WI}\,(p - p_{\text{BHP}}),
\end{equation}
enters the total mass balance~\eqref{eq:total_mass_balance} as a sink term at $\Omega_0^3,$ where $\lambda$ and $p$ denote the density-weighted total mobility and the fluid pressure in $\Omega^3_0$, respectively. The Peaceman well index, $\text{WI}$ is given by
\begin{equation}
 \text{WI} = \frac{2\pi h\, K}{\ln(r_e / r_w) + s},
\end{equation}
where $h$ is the well-cell thickness, $K$ is the permeability, $r_e$ is the Peaceman equivalent radius, and $r_w$ is the wellbore radius, and $s$ is a skin factor (set to zero in this work). The quantity $q_{\text{prod}}$ enters the total mass balance~\eqref{eq:total_mass_balance} as a sink term at $\Omega_0^3$. 

The simulation is carried out over a period of $t_{\text{end}} = 74$ days using an adaptive time-stepping scheme with an initial time-step size $\Delta t_{\mathrm{init}} = 360$ s and a range restricted to $[0.1,\, 360]$ s. At each time step, the nonlinear system \eqref{eq:backward_euler} is solved by Newton's method with the backtracking Armijo line search of Section \ref{sec:armijo_line_search} enabled. The backtracking factor $c$ and Armijo parameter $\theta$ are given in Table \ref{tab:common_parameters}. 

\begin{figure}[H]
    \centering
    \includegraphics[width=0.45\linewidth]{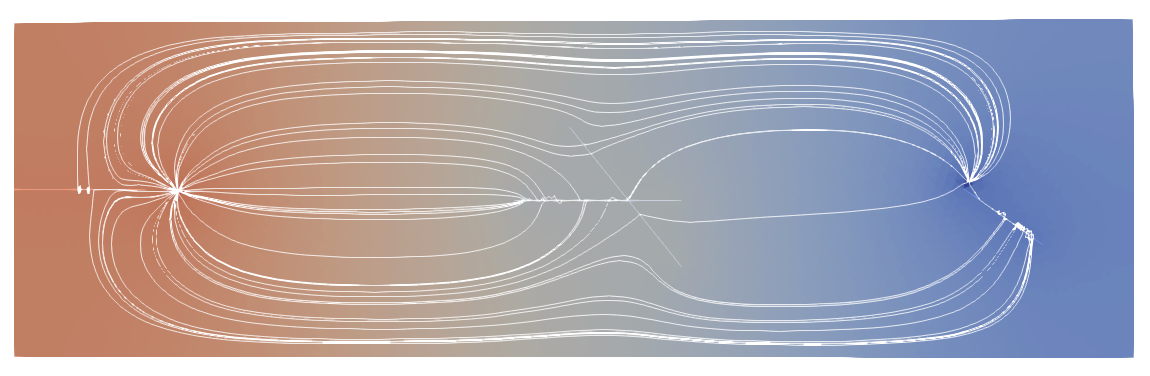}
    \hfill
    \includegraphics[width=0.53\linewidth]{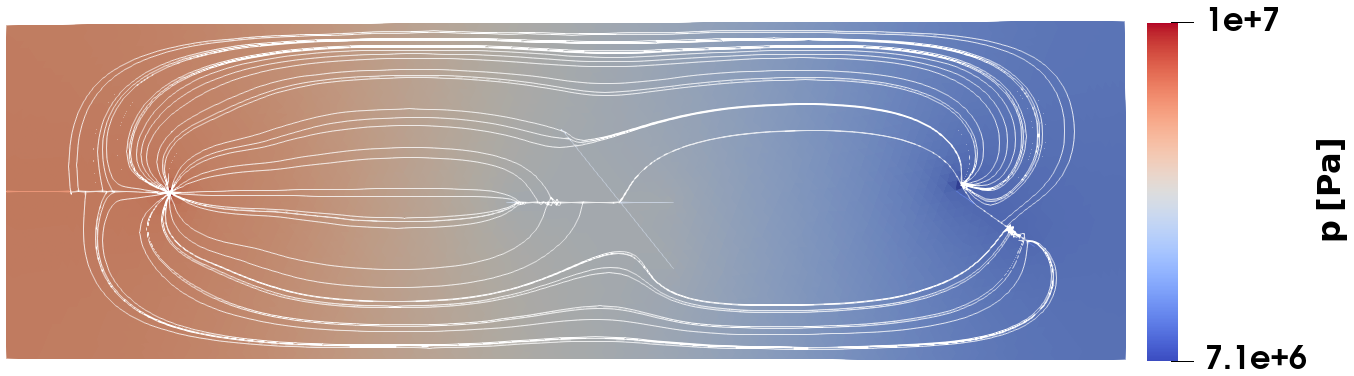}
    \\[0.5em]
    \includegraphics[width=0.45\linewidth]{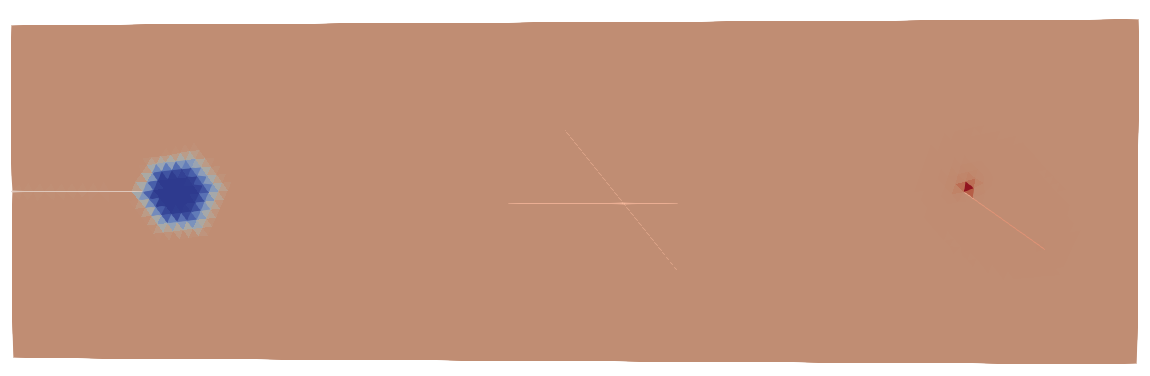}
    \hfill
    \includegraphics[width=0.53\linewidth]{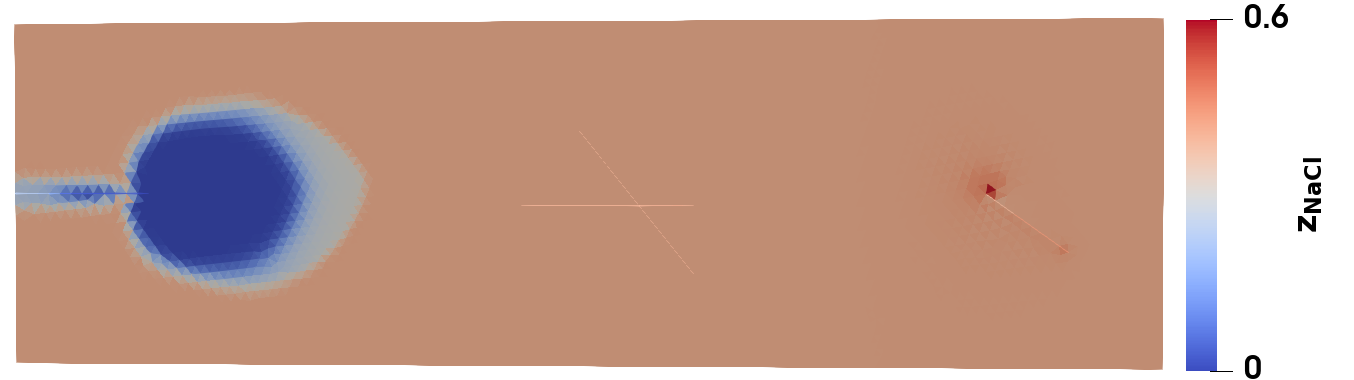}
    \\[0.5em]
    \includegraphics[width=0.45\linewidth]{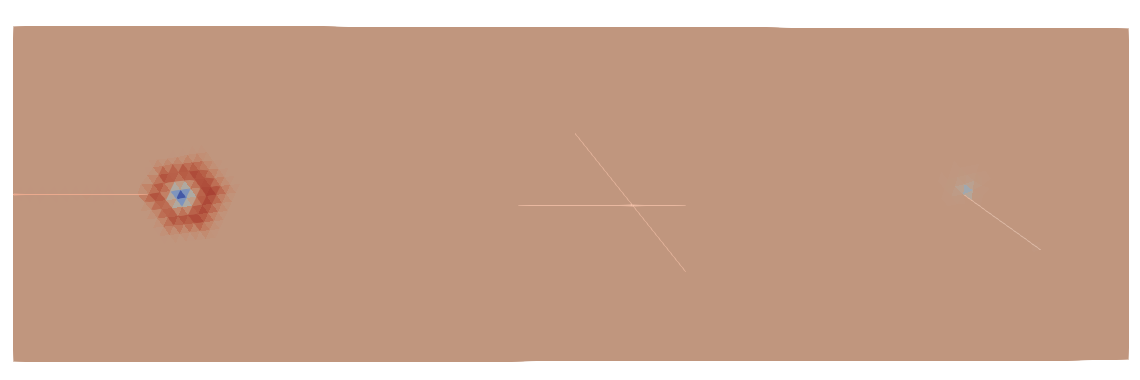}
    \hfill
    \includegraphics[width=0.53\linewidth]{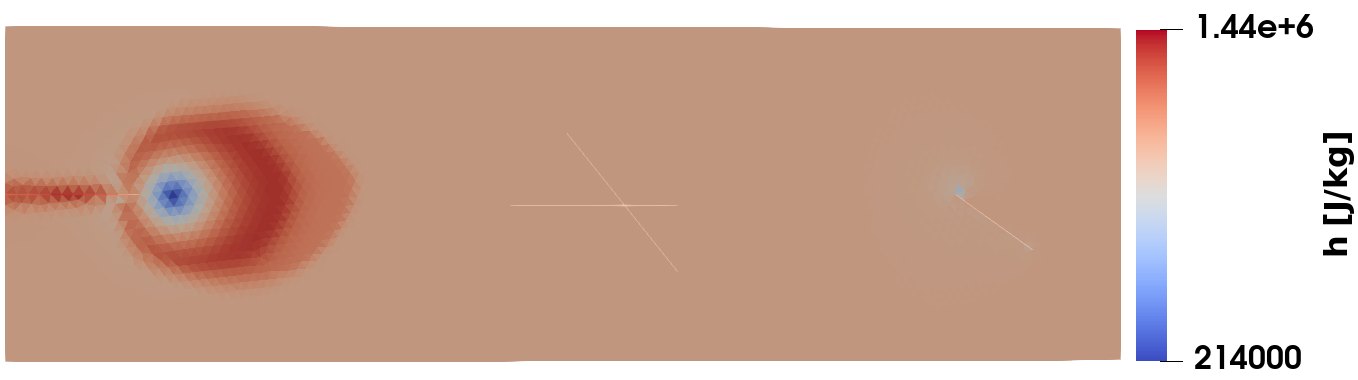}
    \caption{Spatial distribution of the primary variables at $t = 10$\,days (left column) and $t = 74$\,days (right column): pressure $p$ with streamlines (top row), overall salt mass fraction $z_{\mathrm{NaCl}}$ (middle row), and specific enthalpy $h$ (bottom row).}
    \label{fig:primary_variables_profiles}
\end{figure}

\begin{figure}[H]
    \centering
    \includegraphics[width=1.0\linewidth]{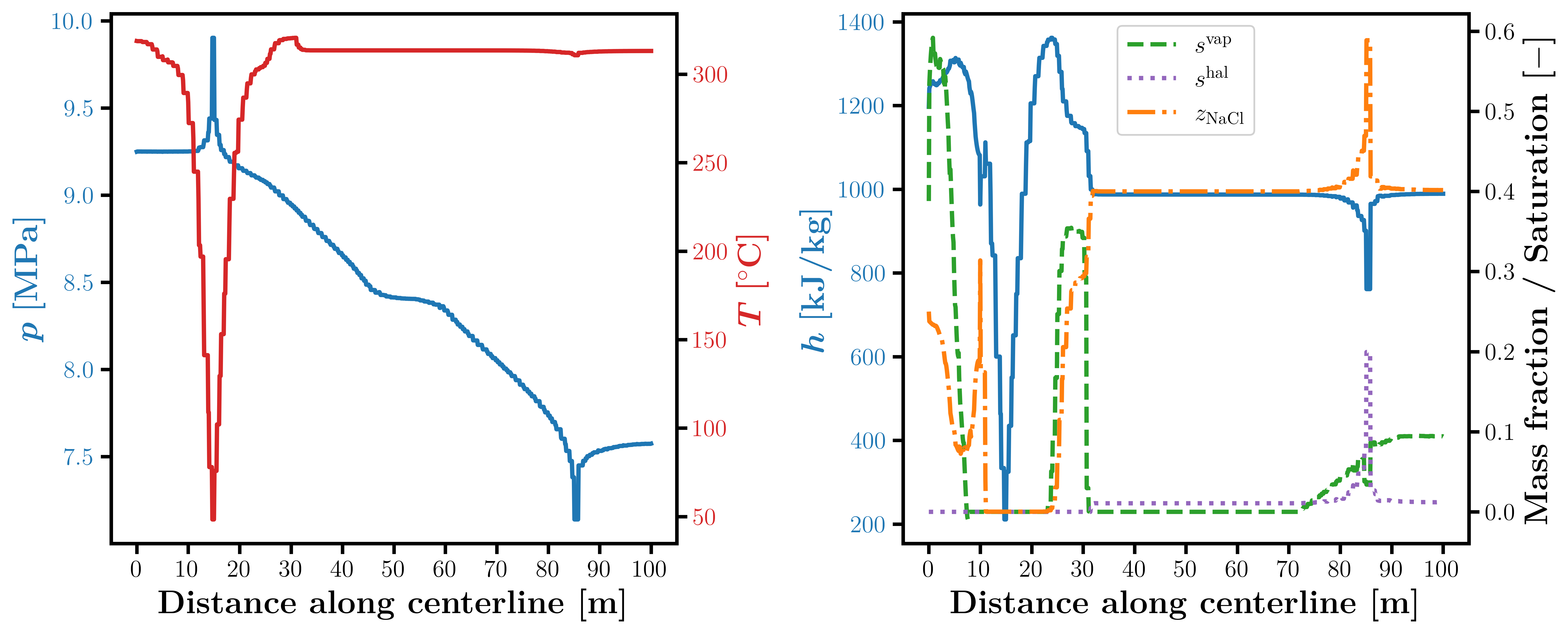}
   \caption{Profiles along the horizontal line $\zeta_2 = 15$ at time $t = 74$\,days. Left: pressure $p$ (blue, left axis) and temperature (red, right axis). Right: specific enthalpy $h$ (solid blue, left axis), alongside vapour saturation $s^{\mathrm{vap}}$ (dashed green), halite saturation $s^{\mathrm{hal}}$ (dotted purple), and overall salt mass fraction $z_\text{NaCl}$ (dash-dotted orange) on the right axis. The injection well is located at $\zeta_1 = 15$\,m and the production well at $\zeta_1 = 85$\,m.}
   \label{fig:central_line_plot_disconnected_0.1}
\end{figure}

Figure~\ref{fig:primary_variables_profiles} displays the spatial distribution of the primary variables at $t = 10$~days (left column) and at the final simulation $t=74$~days (right column). The pressure field (top row) exhibits a smooth gradient through the rock matrix, with localised extrema near the wellbores (see also the pressure profile in the left panel of Figure \ref{fig:central_line_plot_disconnected_0.1}). Because the fracture network is disconnected, the matrix carries the bulk of the inter-well flow. The streamlines, however, show that the discrete fractures locally perturb the flow field by acting as preferential high-permeability conduits. In particular, $\Omega_1^4$ near the producer, $\Omega_1^1$ near the injector, and the intermediate $\Omega_1^2$ and $\Omega_1^3$ capture matrix fluid and channel it into the production region. These results demonstrate the model's ability to capture matrix-fracture flow interactions within mixed-dimensional domain. Comparison between the two columns shows that the pressure field equilibrates rapidly (already in quasi-steady state by $t=10$ days), while the salt and enthalpy fields continue to evolve substantially over the simulation period. 

\begin{figure}[H]
    \centering
    \includegraphics[width=0.48\linewidth]{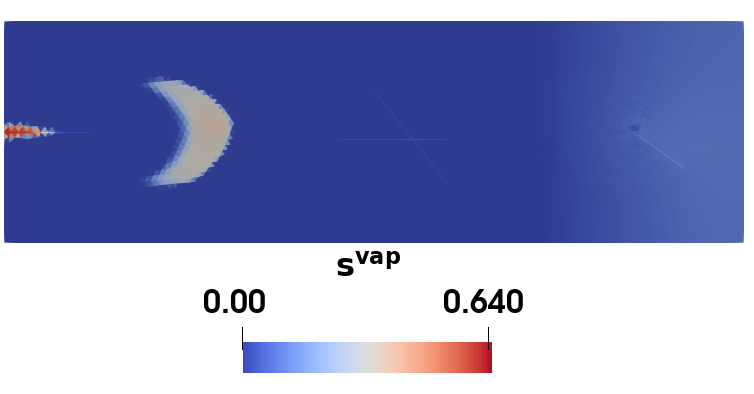}
    \hfill
    \includegraphics[width=0.48\linewidth]{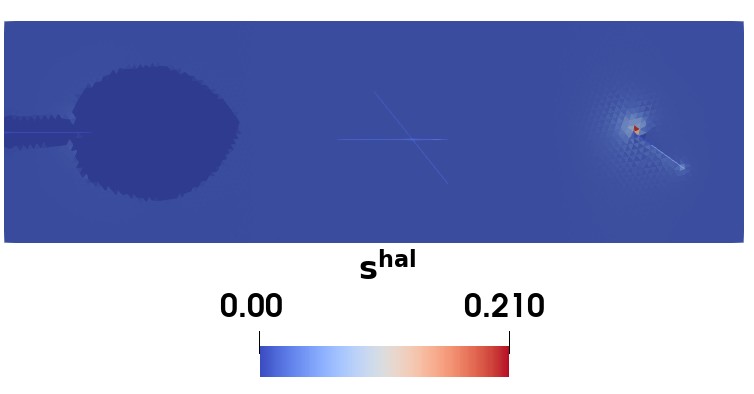}
    \hfill
    \includegraphics[width=0.48\linewidth]{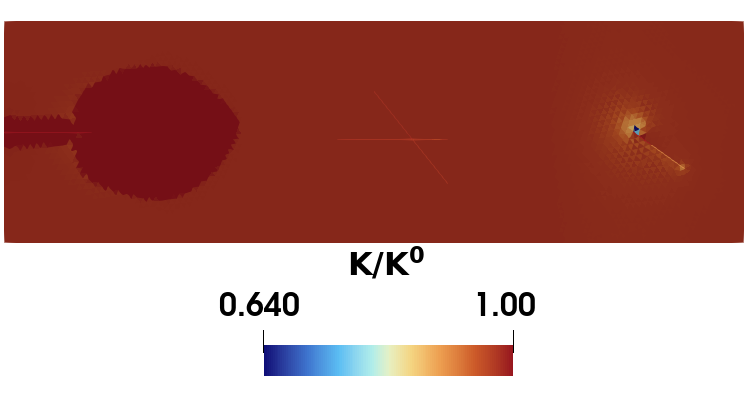}
    \caption{Spatial distribution of secondary fields at $t = 74$\,days: vapour saturation $s^{\mathrm{vap}}$ (top left), halite saturation $s^{\mathrm{hal}}$ (top right), and matrix permeability ratio $K/K^{0}$ (bottom).}
    \label{fig:secondary_fields_ex1}
\end{figure}

The overall salt mass fraction $z_{\text{NaCl}}$ field (Figure~\ref{fig:primary_variables_profiles}, middle row) captures the salt distribution in the reservoir. Around the injector, a pronounced dissolution zone forms, where the in-place solid halite is dissolved in the liquid brine (i.e., $s^\text{hal}=0$, see the right panel of Figure \ref{fig:secondary_fields_ex1}) by the continuous influx of undersaturated liquid water and subsequently advected downstream. The dissolution zone is already established by $t=10$ days as a localised region around the injector and expands substantially by $t=74$ days, developing into two distinct fronts. At the inner front (dark-blue), $z_{\mathrm{NaCl}}$ drops sharply to the injection concentration $z_\text{inj}$ (see $z_{\mathrm{NaCl}}$ profile in the right panel of Figure \ref{fig:central_line_plot_disconnected_0.1}). Dissolution of precipitated halite together with rapid phase redistribution generates strong local gradients in the fluid-mixture enthalpy, producing localised extrema in the enthalpy field that appear as the dark-red ring in Figure~\ref{fig:primary_variables_profiles} (bottom panel). The outer front (light-blue) exhibits different thermodynamic characteristics. Here, the salt fraction remains intermediate ($z_\text{inj}<z_{\text{NaCl}}$), while the local pressure and salinity conditions induce partial boiling of the liquid brine (see Figure \ref{fig:secondary_fields_ex1}, left panel). In this region, the temperature is elevated relative to the cooled inner dissolution zone and coincides with local boiling conditions (see the temperature profile in the left panel of Figure~\ref{fig:central_line_plot_disconnected_0.1}). Similarly, the specific enthalpy differs from that of the inner dissolution zone because of the coupled effects of phase redistribution and latent heat associated with vapour formation (see right panel, Figure~\ref{fig:central_line_plot_disconnected_0.1}).

Conversely, near the production well, a sharp halite precipitation spike appears (Figure \ref{fig:secondary_fields_ex1}, right panel and halite saturation profile in the right panel of Figure \ref{fig:central_line_plot_disconnected_0.1}). Rather than being confined to the wellbore, this precipitation spreads outward to the surrounding rock matrix. The mechanism is the pressure drawdown at the producer, which causes rapid boiling of the liquid brine at the wellbore; boiling concentrates the residual brine until the solubility limit is exceeded, forcing the dissolved salt to precipitate back into solid halite. The precipitation also extends along the adjacent disconnected fracture $\Omega_1^4$ (left panel of Figure~\ref{fig:halite_aperture_ratio_disconnected_0.1}), as the high-permeability conduit propagates the drawdown and boiling zone further into the reservoir. Because halite has a lower specific enthalpy than fluid phases, its accumulation reduces the fluid mixture specific enthalpy, creating the sharp enthalpy dip near the production well (see Figure \ref{fig:primary_variables_profiles}, bottom panel and the enthalpy profile in the right panel of Figure \ref{fig:central_line_plot_disconnected_0.1}).

\begin{figure}[H]
    \centering
    \includegraphics[width=1.0\linewidth]{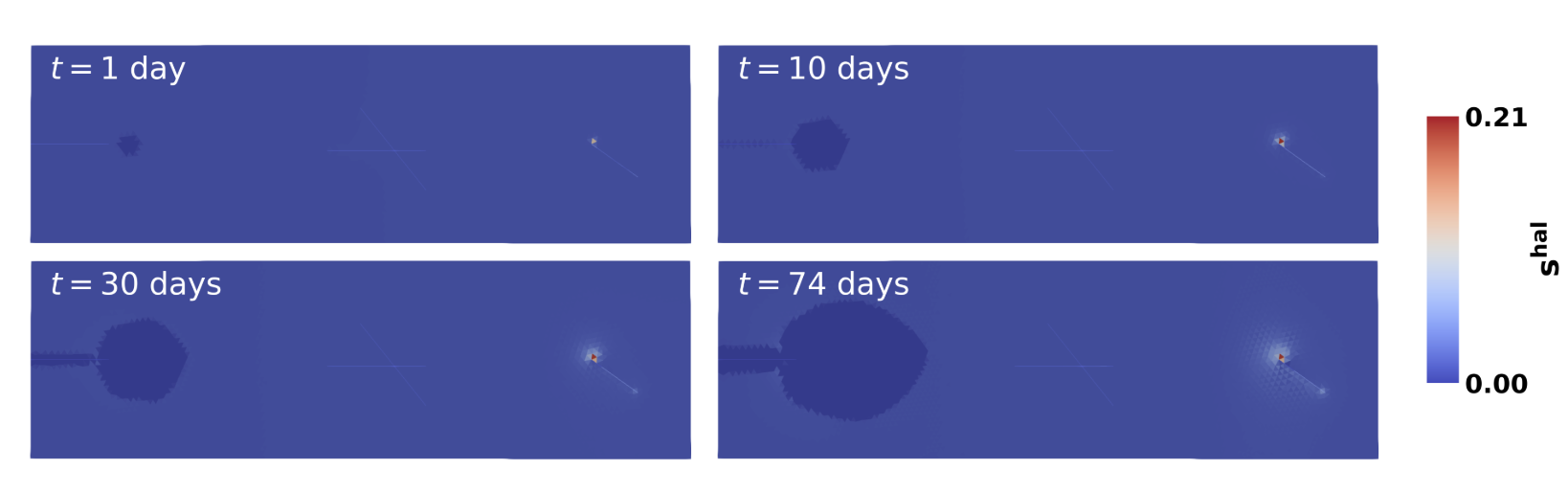}
   \caption{Progression of halite saturation at 1, 10, 30, and 74 days.}
   \label{fig:halite_evolution_profile}
\end{figure}

\begin{figure}[H]
    \centering
    \includegraphics[width=1.0\linewidth]{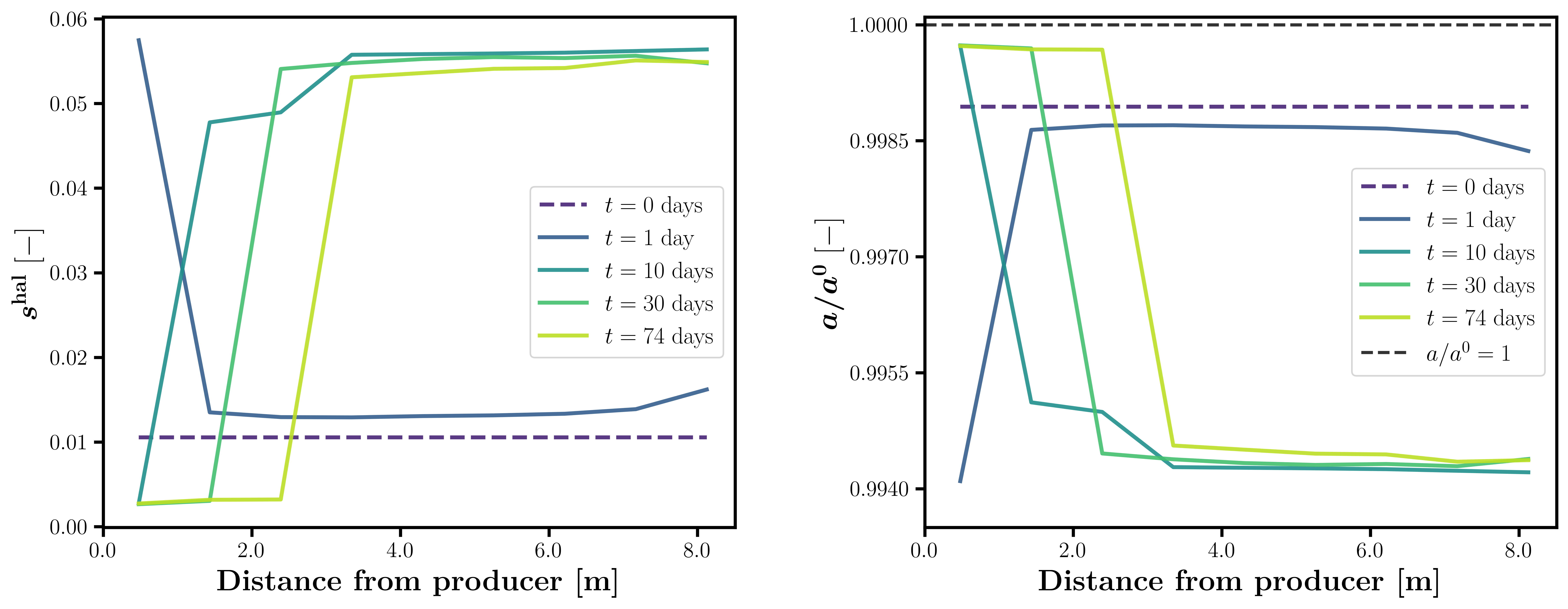}
   \caption{Temporal evolution of halite saturation $s^{\mathrm{hal}}$ (left panel) and aperture ratio $a/a^0$ (right panel) along the fracture $\Omega^4_1$, plotted against distance from the 0D production subdomain $\Omega_0^3$. The producer is at distance $0$; the bottom end of the fracture is at $\approx 8.6\,\text{m}$. Dashed lines indicate the initial state at $t = 0$. In the right panel, the additional black dashed line at $a/a^0 = 1$ marks the hypothetical no-clogging limit.} 
   \label{fig:halite_aperture_ratio_disconnected_0.1}
\end{figure}

The progression of $s^{\mathrm{hal}}$ over time in Figure~\ref{fig:halite_evolution_profile} shows how the two competing zones develop. By $t = 1$\,day, only minor disturbances are visible at the wells. By $t = 10$\,days, the dissolution zone is clearly established at the injector and a precipitation cluster has formed around the producer. The dissolution front continues to advance into the reservoir over time, while the precipitation zone spreads radially from the producer and extends along the adjacent fracture $\Omega^4_1$. The hydraulic consequence of this halite redistribution is shown in Figure~\ref{fig:secondary_fields_ex1} (bottom). Through the
Kozeny--Carman feedback in Equation~\eqref{eq:permeability_kozeny_carman}, the matrix permeability near the producer drops to $K/K^0 \approx 0.64$ in the cells with the highest $s^{\mathrm{hal}}$. Around the injector, $K/K^0$ remains essentially at unity, since the dissolution zone restores the porosity to its halite-free reference value. This near-wellbore permeability reduction also manifests in the streamline pattern shown in Figure~\ref{fig:primary_variables_profiles}. Here, at $t=10$ days, streamlines converge into the producer in a symmetric radial pattern, while at $t=74$ days, a visible deflection appears as the flow re-routes around the spatially non-uniform impaired zone. This illustrates the two-way coupling between flow and precipitation captured by the formulation: the flow drives precipitation through pressure-drawdown-induced boiling, and the resulting permeability reduction in turn modifies the flow field.

To examine the dynamics in the fracture near the production well, we plot $s^{\mathrm{hal}}$ and the aperture ratio $a/a^0$ along $\Omega_1^4$ at several times in Figure~\ref{fig:halite_aperture_ratio_disconnected_0.1}. Initially ($t = 1$~day), rapid boiling causes halite precipitation along the entire fracture, with the highest accumulation in the cell closest to the producer (the throat of $\Omega_1^4$ relative to the production well). As the simulation progresses, however, this precipitation peak moves further down the fracture. In the region near the throat, $s^{\mathrm{hal}}$ gradually drops below the initial value $s_0^{\mathrm{hal}}$, indicating a dissolution front that advances toward the distal end of $\Omega_1^4$. This behaviour reflects the streamline pattern: matrix fluid drawn into $\Omega_1^4$ has been salt-depleted by precipitation in the surrounding matrix, so the fluid arriving at the throat region is undersaturated in salt and dissolves any halite that initially precipitated there. The fracture thus acts as a dissolution conduit, while precipitation occurs primarily in the matrix surrounding the producer and along $\Omega_1^4$. The aperture ratio remains within $0.6\%$ of unity along the entire fracture throughout the simulation, consistent with the mild clogging exponent $\varphi = 0.1$. The sensitivity of this near-fracture dynamics, as well as the reservoir's hydraulic response to the near-well injection rate, is examined in Example~$2$.

\subsection{Example 2}
\label{sec:example2}
This example revisits the disconnected-fracture configuration of Example~1 to assess the sensitivity of the reservoir's hydraulic response to the strength of the aperture--halite feedback. The geometry, mesh, well placement, initial and boundary conditions, and all rock and fluid parameters are retained from Section~\ref{sec:example1}. The only modification is the injection rate, which is increased to ${q}_{\text{inj}} = 0.364$~kgm$^{-3}$s$^{-1}$ (a factor of $1.3$ relative to Example~1). The increased injection rate intensifies near-well boiling and the associated concentration of residual brine, driving stronger halite precipitation around the production well compared to Example 1. Through the aperture-halite feedback of Equations \eqref{eq:aperture_update}-\eqref{eq:cubic_law}, this precipitation amplifies the reduction of fracture aperture and permeability. Here, the final simulation time is $t_{\text{end}} = 7$~days, a duration chosen to capture the onset of fracture-throat clogging.

\begin{figure}[H]
    \centering
    \includegraphics[width=1.0\linewidth]{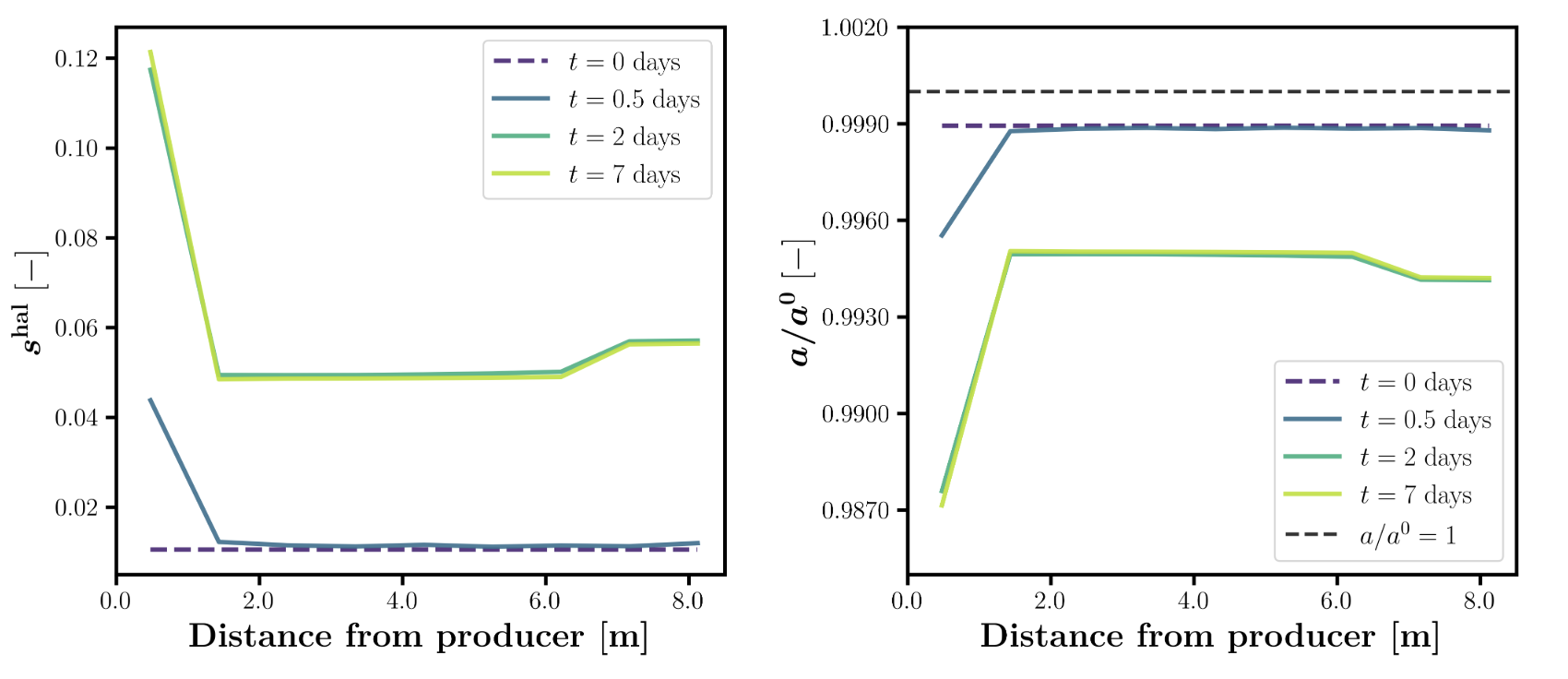}
   \caption{Temporal evolution of halite saturation $s^{\mathrm{hal}}$ (left) and aperture ratio $a/a^0$ (right) along the fracture $\Omega^4_1$ for the ${q}_{\text{inj}} = 0.364$~kgm$^{-3}$s$^{-1}$ case, plotted against distance from the 0D production subdomain $\Omega_0^3$. The producer is at distance $0$; the bottom end of the fracture is at $\approx 8.6\,\text{m}$. Dashed lines indicate the initial state at $t = 0$. In the right panel, the additional black dashed line at $a/a^0 = 1$ marks the hypothetical no-clogging limit. In contrast to Example 1 (Figure \ref{fig:halite_aperture_ratio_disconnected_0.1}), the throat cell now experiences sustained halite precipitation rather than dissolution, with $s^\text{hal}\approx 0.123$ and the aperture ratio dropping to $\approx 0.987$ by $t=7$ days, a larger aperture reduction than the $\approx 0.6\%$ seen at the base injection rate in Example 1.} 
   \label{fig:halite_aperture_ratio_disconnected_0_1_0364}
\end{figure}

The halite saturation and aperture ratio dynamics of the near-production fracture $\Omega_1^4$ are shown in Figure \ref{fig:halite_aperture_ratio_disconnected_0_1_0364}. In contrast to  Example 1, where the fracture throat undergoes dissolution, the higher injection rate now drives the halite precipitation along the entire fracture $\Omega_1^4$. The throat saturation rises over the first two days, from $s^{\mathrm{hal}} \approx 0.04$ at time $t=0.5$ days to a plateau of $s^{\mathrm{hal}} \approx 0.12,$ and remains essentially at that level through $t=7$ days. The increased injection rate thus shifts the fracture from dissolution, observed in Example 1, to sustained precipitation along its full length, driven by intensified near-well boiling that concentrates the residual brine past saturation faster than the undersaturated matrix fluid can redissolve it. The associated aperture reduction nonetheless remains modest: the throat decreases by only about $1.3\%$, compared with the $0.6\%$ of Example 1. Because fracture transmissibility scales with the square of the aperture through the cubic law \eqref{eq:cubic_law}, this small aperture reduction leaves $\Omega_1^4$ essentially open to flow.

\begin{figure}[H]
    \centering
    \includegraphics[width=1.0\linewidth]{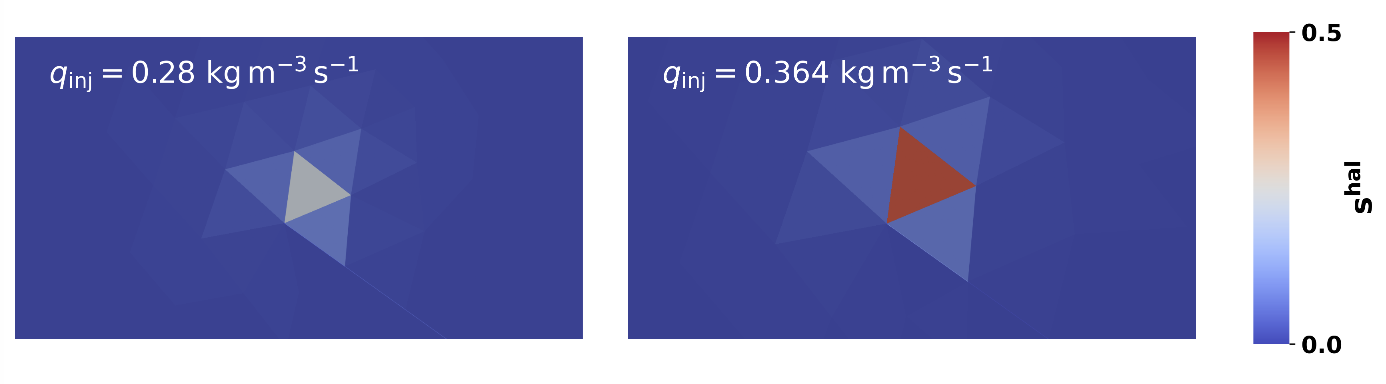}
    \caption{Near-well halite saturation $s^{\text{hal}}$ at $t = 7$~days, zoomed onto the matrix cell hosting the production well and its immediate neighbours. Left: base injection rate $q_{\text{inj}} = 0.28~\text{kg}\,\text{m}^{-3}\,\text{s}^{-1}$. Right: increased rate $q_{\text{inj}} = 0.364~\text{kg}\,\text{m}^{-3}\,\text{s}^{-1}$ ($1.3\times$ base). Halite saturation in the host cell increases from $\approx 0.15$ to $\approx 0.46$.}
   \label{fig:precipitation_disconnected_varphi_1}
\end{figure}

\begin{figure}[H]
    \centering
    \includegraphics[width=1.0\linewidth]{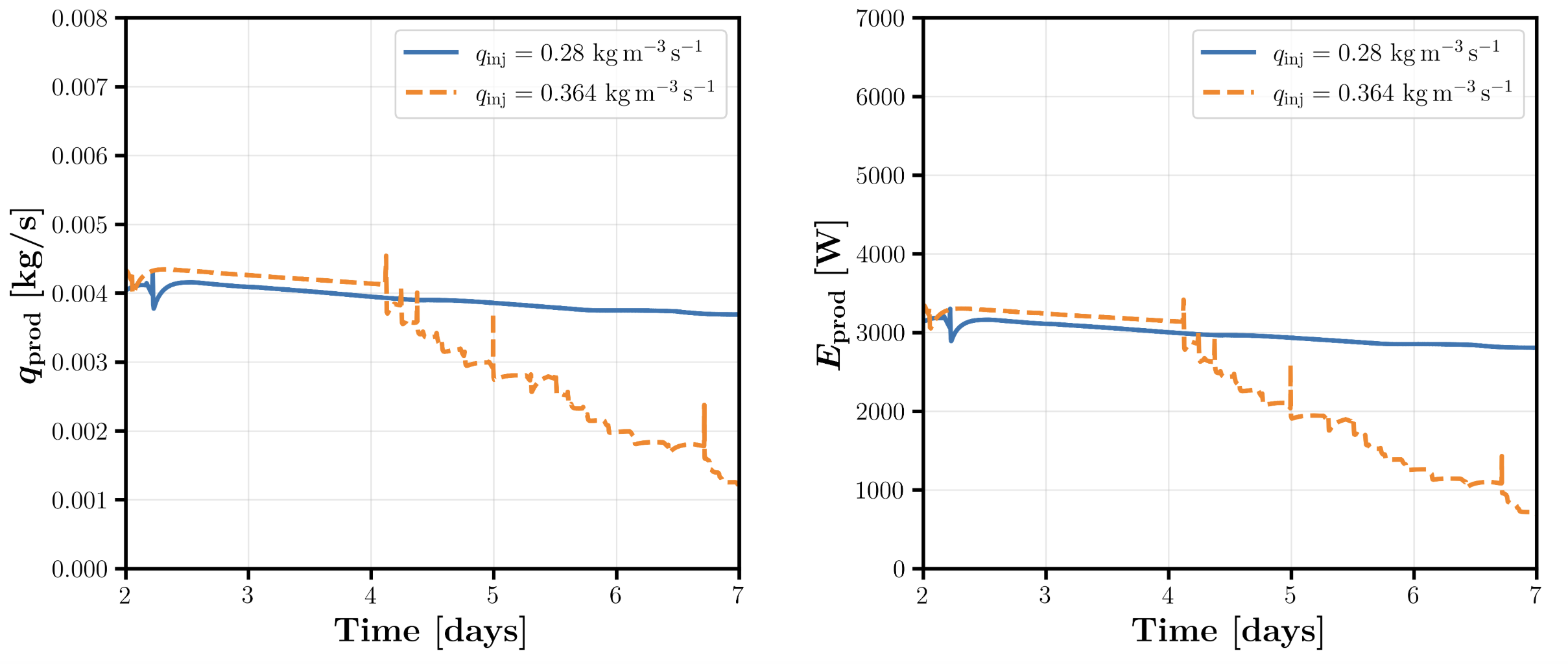}
   \caption{Time evolution of the magnitudes of production mass rate $q_{\text{prod}}$ (left) and energy production rate $E_{\text{prod}}$ (right) for base injection rate (solid) and increased injection rate (dashed), shown for $t \geq 2$~days.} 
   \label{fig:production_diagnostics_disconnected_comparison}
\end{figure}

The increased injection rate has a pronounced consequence on the matrix surrounding the production well. Figure~\ref{fig:precipitation_disconnected_varphi_1} compares the near-wellbore matrix halite saturation at $t = 7$~days for the two injection rates. At the base rate (Example~1), the matrix cell hosting the producer reaches $s^{\mathrm{hal}} \approx 0.15$; at the increased rate, the same cell reaches $s^{\mathrm{hal}} \approx 0.46$. The higher near-well fluid flux intensifies boiling and residual-brine concentration in the host cell, driving this rise in halite saturation. Through the Kozeny--Carman feedback~\eqref{eq:permeability_kozeny_carman}, the host-cell permeability is reduced to $K/K^0 \approx 0.29$, considerably more severe than the $K/K^0 \approx 0.64$ of Example~1. This near-wellbore permeability impairment, rather than the fracture-aperture reduction, produces the substantial decline in both the production mass rate $q_{\text{prod}}$ and the energy production rate $E_{\mathrm{prod}} = q_{\mathrm{prod}} \times h_\mathrm{prod}$, where $h_\mathrm{prod}$ is the specific enthalpy of the produced fluid, shown for the two injection rates in Figure~\ref{fig:production_diagnostics_disconnected_comparison}. 

The clogging exponent $\varphi$ in Equation~\eqref{eq:aperture_update}, which controls the strength of the aperture--halite feedback, has a pronounced effect on the fracture aperture but only a secondary effect on the reservoir-scale response. Increasing $\varphi$ from $0.1$ to $1.0$ at the same injection rate clogs the fracture throat much more strongly, reducing the throat aperture to $a/a^0 \approx 0.88$ (against ${\approx}\,0.987$ at $\varphi = 0.1$) and, through the cubic law~\eqref{eq:cubic_law}, its transmissibility by roughly $23\%$. This partial restriction diverts a larger share of the near-well flow through the adjacent matrix host cell, modestly intensifying boiling there and raising the host-cell halite saturation from $0.46$ to $0.50$ (with $K/K^0$ decreasing from ${\approx}\,0.29$ to ${\approx}\,0.25$).

%
%
%
%

\subsection{Example 3}
\label{sec:example3}
This example examines the role of fracture connectivity by replacing the disconnected fracture network of Examples~1 and 2 with a continuous high-permeability pathway between the wells. The reservoir geometry is illustrated in Figure~\ref{fig:application_geometry} (bottom panel): eight fractures $\Omega_1^1, \ldots, \Omega_1^8$ are arranged in a zigzag chain from the injection well at $(15, 15)$~m to the production well at $(85, 15)$~m, with consecutive segments meeting at seven intersection subdomains $\Omega_0^2, \ldots, \Omega_0^8$. The reference aperture is increased to $a^0 = 10^{-2}$~m, and the injection rate is increased to $q_{\text{inj}} = 0.84~\text{kg}\,\text{m}^{-3}\,\text{s}^{-1}$ (a factor of $3$ relative to Example~1). The clogging exponent is set to $\varphi = 2.0$. All other rock and fluid parameters, well models, and boundary conditions are retained from Section~\ref{sec:example1}. 
The modified parameters are chosen collectively to promote fracture-dominated transport and stronger aperture–precipitation coupling, rather than to isolate the influence of individual parameters. The simulation is advanced to $t_{\text{end}} = 60$~days using an adaptive time-stepping scheme with an initial step $\Delta t^{\text{init}} = 120$~s and a range restricted to $[1, 600]$~s. 

\begin{figure}[H]
    \centering
    \includegraphics[width=0.449\linewidth]{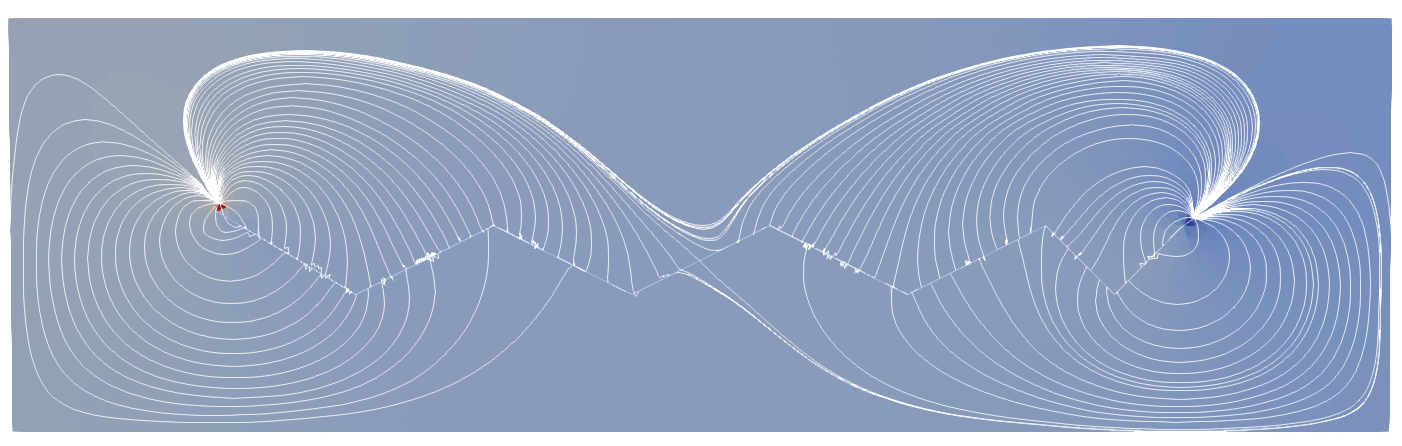}
    \hfill
    \includegraphics[width=0.53\linewidth]{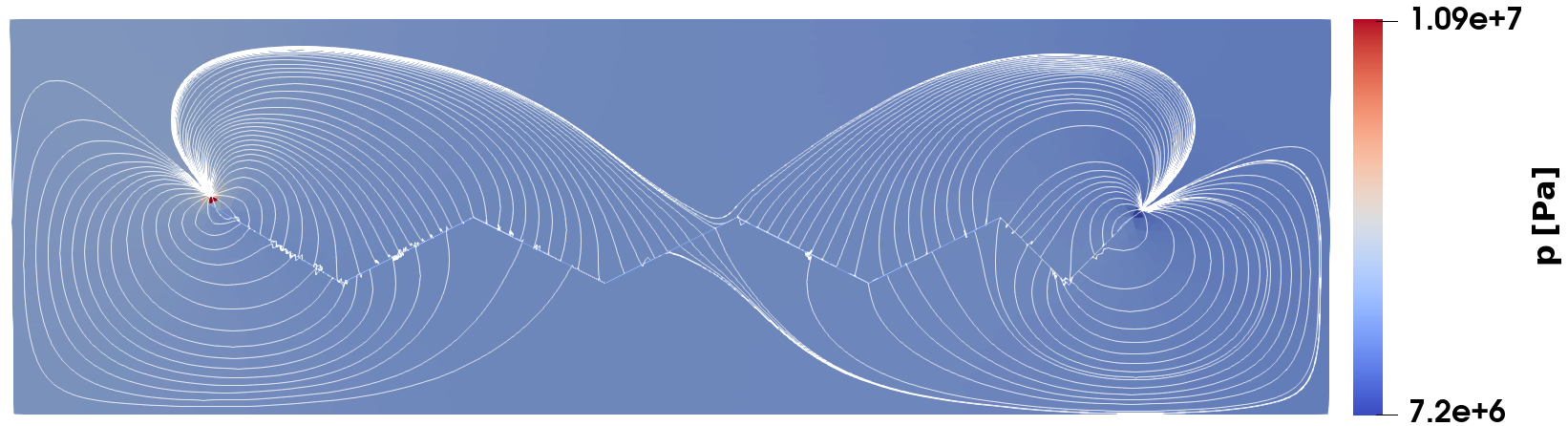}
    \\[0.5em]
    \includegraphics[width=0.449\linewidth]{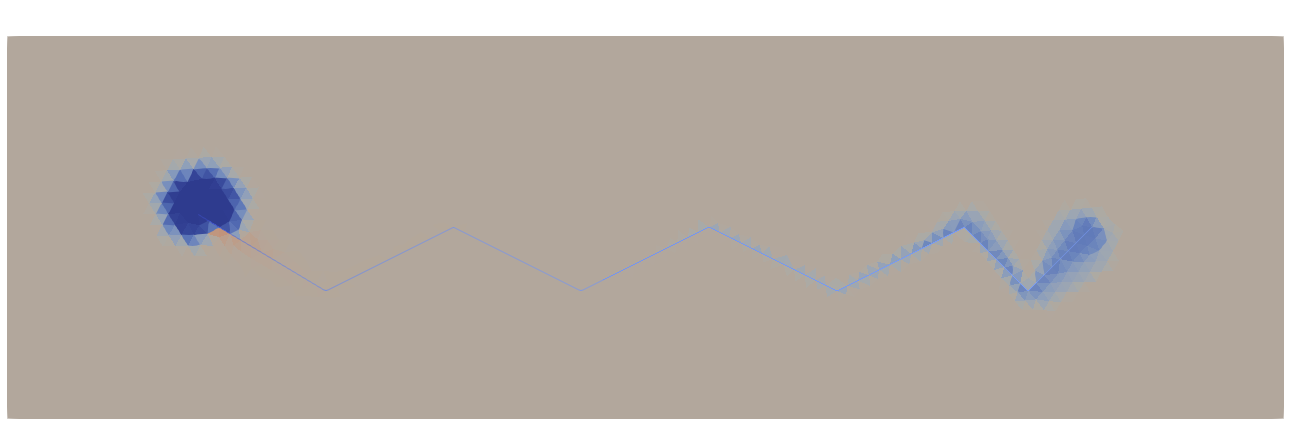}
    \hfill
    \includegraphics[width=0.53\linewidth]{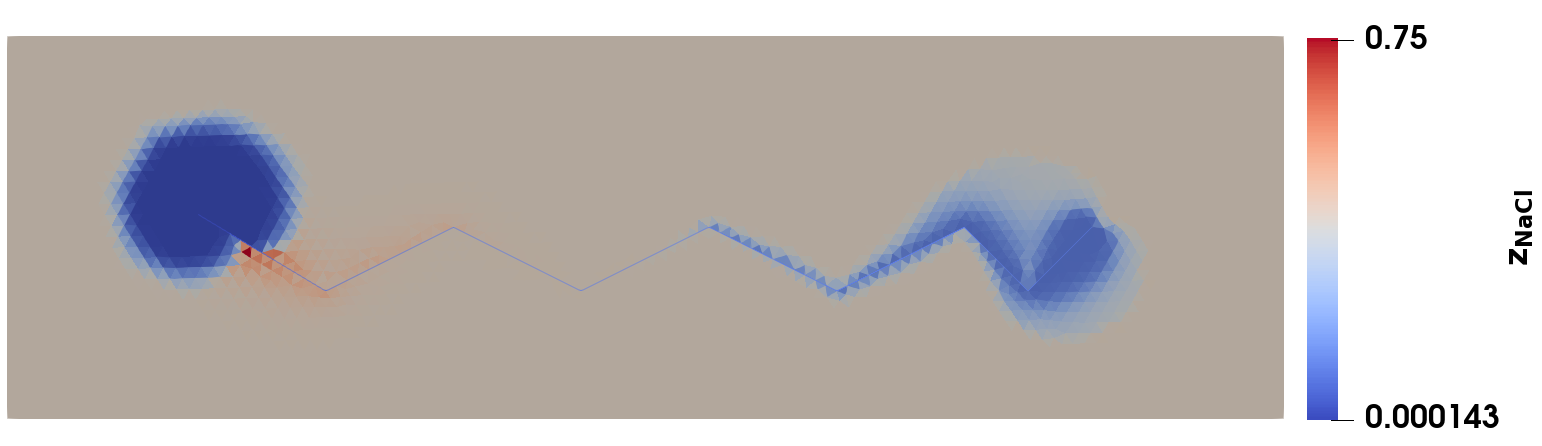}
    \\[0.5em]
    \includegraphics[width=0.449\linewidth]{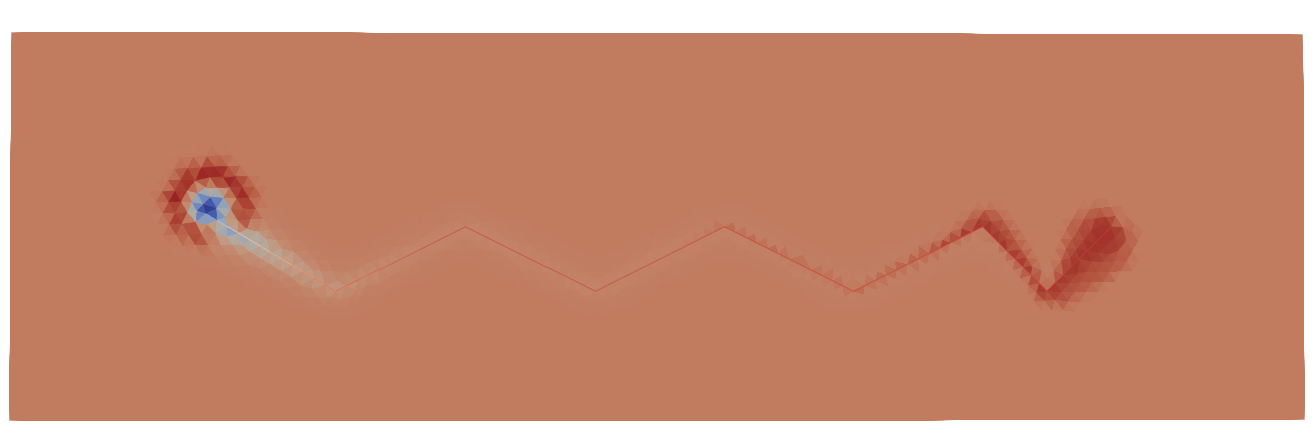}
    \hfill
    \includegraphics[width=0.53\linewidth]{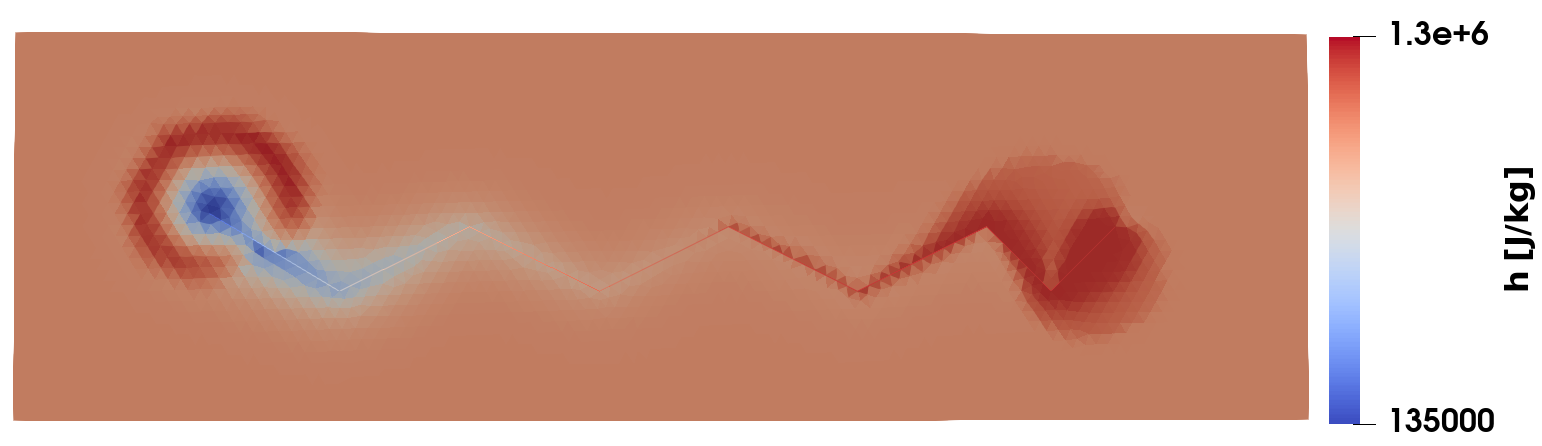}
    \caption{Spatial distribution of the primary variables at $t = 10$\,days (left column) and $t = 60$\,days (right column): pressure $p$ with streamlines (top row), overall salt mass fraction $z_{\mathrm{NaCl}}$ (middle row), and specific enthalpy $h$ (bottom row).}
    \label{fig:primary_variables_profiles_connected}
\end{figure}

\begin{figure}[H]
    \centering
    \includegraphics[width=1.0\linewidth]{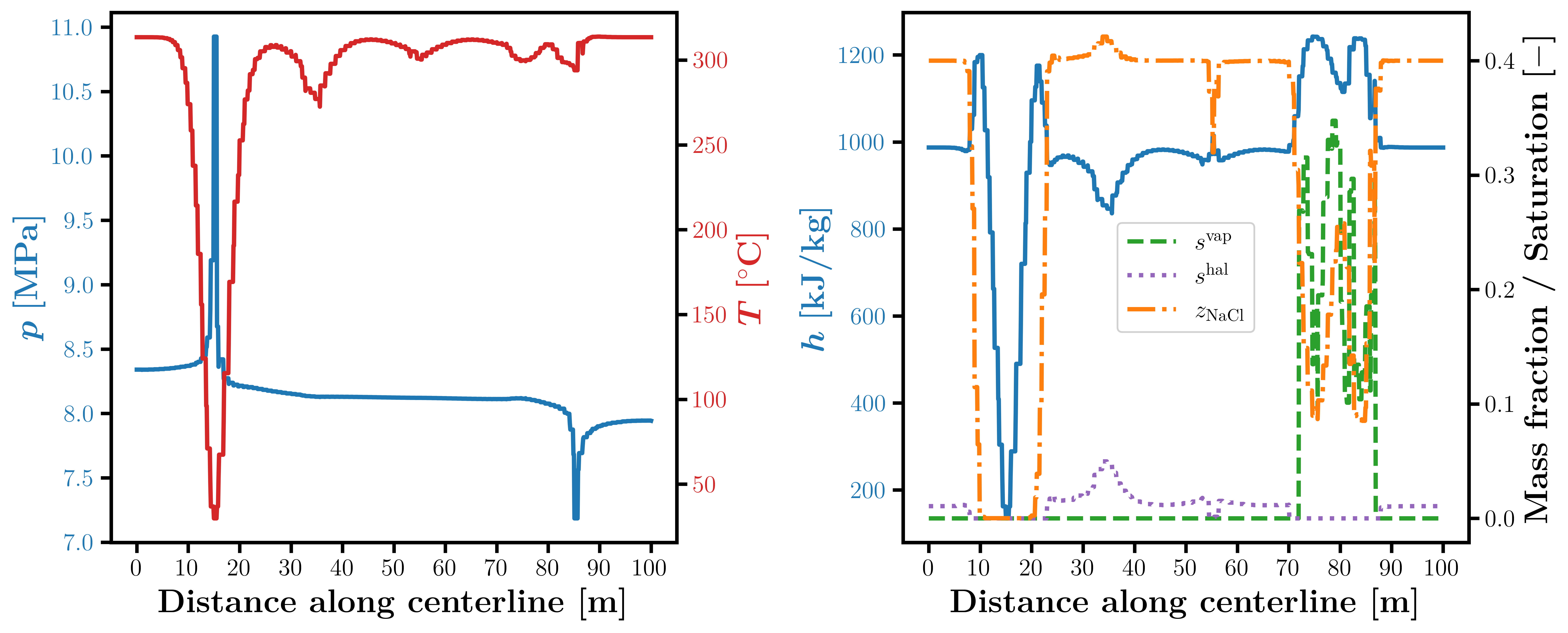}
   \caption{Profiles along the horizontal line $\zeta_2 = 16$ at time $t = 60$\,days. Left: pressure $p$ (blue, left axis) and temperature (red, right axis). Right: specific enthalpy $h$ (solid blue, left axis), alongside vapour saturation $s^{\mathrm{vap}}$ (dashed green), halite saturation $s^{\mathrm{hal}}$ (dotted purple), and overall salt mass fraction $z_\text{NaCl}$ (dash-dotted orange) on the right axis. The injection well is located at $\zeta_1 = 15$\,m and the production well at $\zeta_1 = 85$\,m. The oscillatory pattern between $\zeta_1 \approx 25$ and $70$~m is a sampling artefact: the line $\zeta_2 = 16$ cuts through the zigzag fracture chain, producing local dips in the profiles near the successive intersections.}
   \label{fig:central_line_plot_connected_2.0}
\end{figure}

Figure~\ref{fig:primary_variables_profiles_connected} displays the spatial distribution of the primary variables at $t=10$~days (left column) and at the final simulation time of $60$~days (right column). The pressure field (top row) retains similar local extrema observed in the disconnected cases (see also the pressure profile in the left panel of Figure~\ref{fig:central_line_plot_connected_2.0}). However, between the two wells, the pressure is approximately uniform, because the high conductivity of the fracture chain allows fluid flow with only a small pressure gradient along the fracture chain. The streamlines show that inter-well transport is dominated by the fracture chain. The pressure field is essentially identical between two columns, reflecting rapid pressure equilibration through the high-conductivity fracture pathway.

The overall salt mass fraction field (Figure \ref{fig:primary_variables_profiles_connected}, middle row) shows a region around the injector that is salt-depleted due to the continuous influx of low-salinity water displacing the in-place halite-saturated brine (see also the $z_{\text{NaCl}}$ profile in the right panel of Figure \ref{fig:central_line_plot_connected_2.0}). In the surrounding matrix along the upstream fractures $\Omega_1^1, \Omega_1^2$, the salinity is elevated, reflecting the advective transport of salt displaced from upstream. Around the production well, a second salt-depleted region develops, where the low-salinity fluid arriving through the fracture chain mixes with the resident brine. The fracture chain itself carries low-salinity fluid throughout the simulation, maintaining $z_{\text{NaCl}}$ below the initial value along its full length. Both salt-depleted regions grow substantially between $t=10$ and $t = 60$ days, with the producer-side depletion expanding as the arriving low-salinity fluid progressively displaces the resident brine.

\begin{figure}[H]
    \centering
    \includegraphics[width=0.48\linewidth]{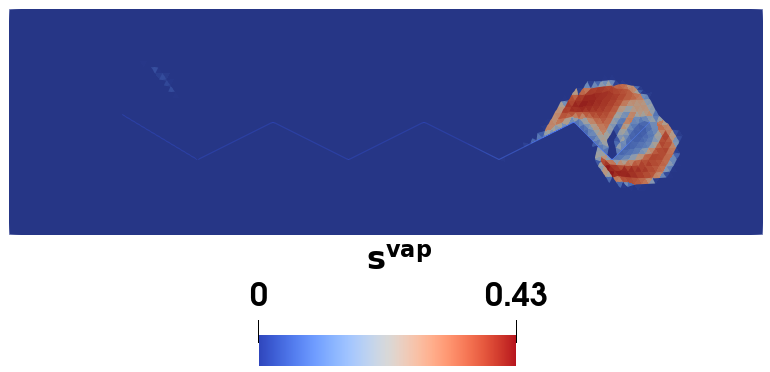}
    \hfill
    \includegraphics[width=0.48\linewidth]{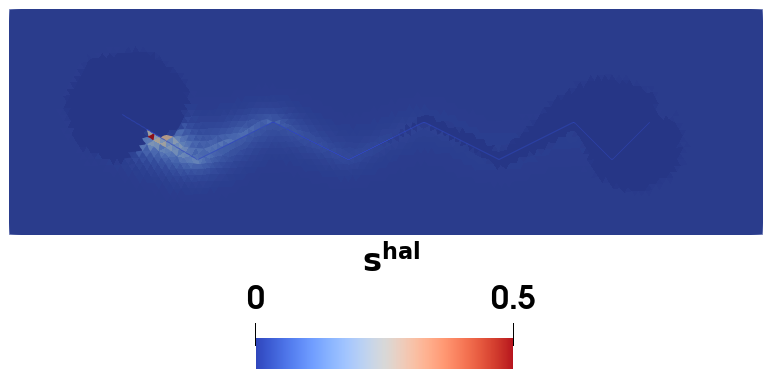}
    \caption{Spatial distribution of the phase saturations at $t = 60$\,days: vapour saturation $s^{\mathrm{vap}}$ (left) and halite saturation $s^{\mathrm{hal}}$ (right).}
    \label{fig:halite_vapor_saturation_profile_connected}
\end{figure}

The halite saturation field (Figure~\ref{fig:halite_vapor_saturation_profile_connected}, right panel) shows that the in-place solid halite has fully dissolved both around the injector and around the producer (see also the $s^{\text{hal}}$ profile in the right panel of Figure \ref{fig:central_line_plot_connected_2.0}). The dissolution around the producer is driven by the low-salinity fluid delivered through the fracture chain, which consumes the in-place halite. Halite precipitation occurs in two distinct regions: a thin band immediately ahead of the injector-side dissolution, where the matrix-resident brine is enriched by advected NaCl from the depleted zone and cooled by contact with the colder upstream fluid, exceeding the local halite solubility; and a streak along the fracture chain, where the cool low-salinity fluid in the chain cools the adjacent matrix through fracture--matrix heat exchange, reducing the local halite solubility of the matrix-resident brine and driving precipitation in the matrix along the chain. 

The vapour saturation field (Figure~\ref{fig:halite_vapor_saturation_profile_connected}, left panel) shows that boiling is concentrated in the matrix cells surrounding the producer. Localised vapour formation is also visible along the producer-side fracture segments, reflecting the combined effects of depressurisation and progressive thermal equilibration of the flowing brine with the surrounding hot matrix. Despite the intense boiling around the producer, no halite precipitates there. Although boiling enriches the residual liquid brine in salt, the fluid transported through the fracture chain remains sufficiently dilute that the local brine salinity does not reach halite saturation. This contrasts with the boiling-induced precipitation observed near the producer in Examples~1 and~2.

The specific enthalpy distribution is shown in Figure~\ref{fig:primary_variables_profiles_connected} (bottom row). Near the injector, low-enthalpy injected fluid creates a localised cooling region surrounded by sharp enthalpy gradients associated with dissolution and phase redistribution. Along the connected fracture chain, the specific enthalpy increases progressively from injector to producer as the flowing brine absorbs heat from the surrounding hot rock matrix. Near the producer, the enthalpy exceeds the initial reservoir value because depressurisation and phase separation generate a vapour-saturated, low-salinity fluid mixture with elevated specific enthalpy (see the enthalpy profile in the right panel of Figure \ref{fig:central_line_plot_connected_2.0}). This contrasts with the enthalpy drop near the producer in Example~1.

\section{Conclusion}
\label{sec:conclusion}

We develop a new compositional framework for non-isothermal multiphase flow with salt precipitation and dissolution in high-enthalpy fractured geothermal reservoirs. The formulation combines a persistent-variable approach, which handles phase appearance and disappearance within a fixed set of primary variables, with a discrete fracture-matrix (DFM) representation of the fractured reservoir. The H$_2$O--NaCl thermodynamics follow the Driesner correlations and are efficiently evaluated through operator-based multilinear interpolation. Porosity, fracture aperture, and permeability evolve with dissolution and precipitation through Kozeny--Carman-type feedback relations. The governing mixed-dimensional equation system is discretised in space by a cell-centered finite volume scheme with multi-point flux approximation for the diffusive fluxes and upwinding for the advective terms, and fully implicitly in time using backward Euler with adaptive step control. The resulting nonlinear system is solved using Newton's method with a backtracking Armijo line search.

The model is implemented in the open-source PorePy framework and verified against the CSMP++ simulator on a one-dimensional salt-dissolution benchmark spanning four phase regions, including vapour+halite, vapour+liquid, single-phase liquid, and liquid+halite. Despite the strongly nonlinear thermodynamic properties, the two simulators, which employ different numerical approaches, show excellent agreement across all phase regions and their transitions.

The production-driven numerical examples demonstrate that fracture connectivity significantly impacts the thermo--hydraulic--compositional feedbacks in high-enthalpy geothermal reservoirs. In the disconnected fracture configurations, inter-well transport occurs primarily through the matrix, promoting near-production well salt precipitation and permeability impairment. Increasing the injection rate intensifies localised fracture-throat clogging and substantially reduces production performance. In contrast, the connected fracture network promotes rapid fracture-dominated transport and stronger pressure communication between the wells, while simultaneously redistributing precipitation away from the production region. Although boiling remains pronounced near the production well, the transported brine remains sufficiently dilute so that halite saturation is not reached locally. These results show that fracture connectivity controls not only fluid transport pathways, but also the spatial distribution of precipitation-driven transmissibility reduction.


\section*{Data availability}

The data and source code for the results presented in this work are available,
and the simulations and figures can be reproduced using a Docker container
hosted at: \url{https://doi.org/10.5281/zenodo.20578448}.

\section*{Acknowledgments}

This project has received funding from the European Research Council (ERC)
under the European Union's Horizon 2020 research and innovation programme
(grant agreement No.\ 101002507).

\bibliographystyle{cas-model2-names}
\bibliography{micRef}

@techreport{IEA2023,
  author = {{International Energy Agency}},
  title = {World Energy Outlook 2023},
  institution = {International Energy Agency},
  year = {2023},
  address = {Paris, France},
  url = {https://www.iea.org/reports/world-energy-outlook-2023},
  note = {Accessed: 2026-01-09}
}

@article{xing2017parallel,
  title={Parallel numerical modeling of hybrid-dimensional compositional non-isothermal Darcy flows in fractured porous media},
  author={Xing, Feng and Masson, Roland and Lopez, Simon},
  journal={Journal of Computational Physics},
  volume={345},
  pages={637--664},
  year={2017},
  publisher={Elsevier}
}

@article{aghili2021hybrid,
  title={A hybrid-dimensional compositional two-phase flow model in fractured porous media with phase transitions and Fickian diffusion},
  author={Aghili, Joubine and De Dreuzy, Jean-Raynald and Masson, Roland and Trenty, Laurent},
  journal={Journal of Computational Physics},
  volume={441},
  pages={110452},
  year={2021},
  publisher={Elsevier}
}

@article{boon2018robust,
  title={Robust discretization of flow in fractured porous media},
  author={Boon, Wietse M and Nordbotten, Jan M and Yotov, Ivan},
  journal={SIAM Journal on Numerical Analysis},
  volume={56},
  number={4},
  pages={2203--2233},
  year={2018},
  publisher={SIAM}
}

@article{dugstad2022dimensional,
  title={Dimensional reduction of a fractured medium for a two-phase flow},
  author={Dugstad, Martin and Kumar, Kundan},
  journal={Advances in Water Resources},
  volume={162},
  pages={104140},
  year={2022},
  publisher={Elsevier}
}

@article{maimoni1982minerals,
  title={Minerals recovery from Salton Sea geothermal brines: A literature review and proposed cementation process},
  author={Maimoni, A},
  journal={Geothermics},
  volume={11},
  number={4},
  pages={239--258},
  year={1982},
  publisher={Elsevier}
}

@article{grant2020trace,
  title={Trace metal distributions in sulfide scales of the seawater-dominated Reykjanes geothermal system: Constraints on sub-seafloor hydrothermal mineralizing processes and metal fluxes},
  author={Grant, Hannah LJ and Hannington, Mark D and Hardard{\'o}ttir, Vigd{\'\i}s and Fuchs, Sebastian H and Schumann, Dirk},
  journal={Ore Geology Reviews},
  volume={116},
  pages={103145},
  year={2020},
  publisher={Elsevier}
}

@article{ingridgeothermal,
  title={Geothermal Energy-From Theoretical Models to Exploration and Development.},
  author={Stober Ingrid, KB},
  year={2021},
  journal={Switzerland: Springer Cham}
}

@misc{oguntola2026docker,
  author       = {Oguntola, Micheal B. and 
                  Duran, Omar and 
                  Keilegavlen, Eirik and 
                  Berre, Inga},
  title        = {Source code: Mathematical Modeling of Salt Precipitation and Multi-Phase Flow in High Enthalpy Fractured Geothermal Systems (v1.0.0)},
  year         = {2026},
  publisher    = {Zenodo},
  doi          = {10.5281/zenodo.20451960},
  url          = {https://doi.org/10.5281/zenodo.20578448}
}

@article{hesshaus2013halite,
  title={Halite clogging in a deep geothermal well--Geochemical and isotopic characterisation of salt origin},
  author={Hesshaus, Annalena and Houben, Georg and Kringel, Robert},
  journal={Physics and Chemistry of the Earth, Parts A/B/C},
  volume={64},
  pages={127--139},
  year={2013},
  publisher={Elsevier}
}

@book{phillips1991flow,
  title={Flow and reactions in permeable rocks},
  author={Phillips, Owen M},
  year={1991},
  publisher={Cambridge University Press}
}

@article{gunnlaugsson2012scaling,
  title={Scaling in geothermal installation in Iceland},
  author={Gunnlaugsson, Einar},
  journal={Proceedings of short course on geothermal development and geothermal wells. Santa Tecla, El Salvado},
  year={2012}
}

@article{scott2017boiling,
  title={Boiling and condensation of saline geothermal fluids above magmatic intrusions},
  author={Scott, Samuel and Driesner, Thomas and Weis, Philipp},
  journal={Geophysical Research Letters},
  volume={44},
  number={4},
  pages={1696--1705},
  year={2017},
  publisher={Wiley Online Library}
}

@article{cavarretta1990schorl,
  title={Schorl-dravite-ferridravite tourmalines deposited by hydrothermal magmatic fluids during early evolution of the Larderello geothermal field, Italy},
  author={Cavarretta, Giuseppe and Puxeddu, Mariano},
  journal={Economic Geology},
  volume={85},
  number={6},
  pages={1236--1251},
  year={1990},
  publisher={Society of Economic Geologists}
}

@incollection{von2016silica,
  title={Silica scale control in geothermal plants—Historical perspective and current technology},
  author={von Hirtz, Paul},
  booktitle={Geothermal power generation},
  pages={443--476},
  year={2016},
  publisher={Elsevier}
}

@article{hardardottir2010cu,
  title={Cu-rich scales in the Reykjanes geothermal system, Iceland},
  author={Hardard{\'o}ttir, Vigd{\'\i}s and Hannington, Mark and Hedenquist, Jeffrey and Kjarsgaard, Ingrid and Hoal, Karin},
  journal={Economic Geology},
  volume={105},
  number={6},
  pages={1143--1155},
  year={2010},
  publisher={Society of Economic Geologists}
}

@article{ji2025capillary,
  title={Capillary-Driven transport and precipitation of salt in heterogeneous structures during carbon sequestration},
  author={Ji, Tiancheng and Haghi, Amir H and Jiang, Peixue and Chalaturnyk, Rick and Xu, Ruina},
  journal={Geophysical Research Letters},
  volume={52},
  number={13},
  pages={e2024GL114388},
  year={2025},
  publisher={Wiley Online Library}
}

@article{afanasyev2018formation,
  title={Formation of magmatic brine lenses via focussed fluid-flow beneath volcanoes},
  author={Afanasyev, Andrey and Blundy, Jon and Melnik, Oleg and Sparks, Steve},
  journal={Earth and Planetary Science Letters},
  volume={486},
  pages={119--128},
  year={2018},
  publisher={Elsevier}
}

@article{samardzioska2005numerical,
  title={Numerical comparison of the equivalent continuum, non-homogeneous and dual porosity models for flow and transport in fractured porous media},
  author={Samardzioska, Todorka and Popov, Viktor},
  journal={Advances in water resources},
  volume={28},
  number={3},
  pages={235--255},
  year={2005},
  publisher={Elsevier}
}

@article{les2025geothermal,
  title={Geothermal modeling in complex geological systems with ComPASS},
  author={Les Landes, A Armandine and Beaude, Laurence and Quiroz, D Castanon and Jeannin, Laurent and Lopez, Simon and Sma{\"\i}, Farid and Guillon, Th{\'e}ophile and Masson, Roland},
  journal={Computers \& Geosciences},
  volume={194},
  pages={105752},
  year={2025},
  publisher={Elsevier}
}

@article{li2019coupled,
  title={Coupled thermo-hydro-mechanical analysis of stimulation and production for fractured geothermal reservoirs},
  author={Li, Sanbai and Feng, Xia-Ting and Zhang, Dongxiao and Tang, Huiying},
  journal={Applied Energy},
  volume={247},
  pages={40--59},
  year={2019},
  publisher={Elsevier}
}

@article{zeng2013numerical,
  title={Numerical simulation of heat production potential from hot dry rock by water circulating through a novel single vertical fracture at Desert Peak geothermal field},
  author={Zeng, Yu-Chao and Wu, Neng-You and Su, Zheng and Wang, Xiao-Xing and Hu, Jian},
  journal={Energy},
  volume={63},
  pages={268--282},
  year={2013},
  publisher={Elsevier}
}

@article{martin2005modeling,
  title={Modeling fractures and barriers as interfaces for flow in porous media},
  author={Martin, Vincent and Jaffr{\'e}, J{\'e}r{\^o}me and Roberts, Jean E},
  journal={SIAM Journal on Scientific Computing},
  volume={26},
  number={5},
  pages={1667--1691},
  year={2005},
  publisher={SIAM}
}

@article{jiang2015multimechanistic,
  title={A multimechanistic multicontinuum model for simulating shale gas reservoir with complex fractured system},
  author={Jiang, Jiamin and Younis, Rami M},
  journal={Fuel},
  volume={161},
  pages={333--344},
  year={2015},
  publisher={Elsevier}
}

@article{berre2019flow,
  title={Flow in fractured porous media: A review of conceptual models and discretization approaches},
  author={Berre, Inga and Doster, Florian and Keilegavlen, Eirik},
  journal={Transport in Porous Media},
  volume={130},
  number={1},
  pages={215--236},
  year={2019},
  publisher={Springer}
}

@article{gao2022review,
  title={A review of simulation models of heat extraction for a geothermal reservoir in an enhanced geothermal system},
  author={Gao, Xiang and Li, Tailu and Zhang, Yao and Kong, Xiangfei and Meng, Nan},
  journal={Energies},
  volume={15},
  number={19},
  pages={7148},
  year={2022},
  publisher={MDPI}
}

@article{chen2024capillary,
  title={Capillary-driven backflow during salt precipitation in a rough fracture},
  author={Chen, Xu-Sheng and Hu, Ran and Zhou, Chen-Xing and Xiao, Yang and Yang, Zhibing and Chen, Yi-Feng},
  journal={Water Resources Research},
  volume={60},
  number={3},
  pages={e2023WR035451},
  year={2024},
  publisher={Wiley Online Library}
}

@article{noiriel2021geometry,
  title={Geometry and mineral heterogeneity controls on precipitation in fractures: An X-ray micro-tomography and reactive transport modeling study},
  author={Noiriel, Catherine and Seigneur, Nicolas and Le Guern, Pierre and Lagneau, Vincent},
  journal={Advances in Water Resources},
  volume={152},
  pages={103916},
  year={2021},
  publisher={Elsevier}
}

@article{nooraiepour2018effect,
  title={Effect of CO2 phase states and flow rate on salt precipitation in shale caprocks—a microfluidic study},
  author={Nooraiepour, Mohammad and Fazeli, Hossein and Miri, Rohaldin and Hellevang, Helge},
  journal={Environmental science \& technology},
  volume={52},
  number={10},
  pages={6050--6060},
  year={2018},
  publisher={ACS Publications}
}

@article{hyman2020flow,
  title={Flow channeling in fracture networks: characterizing the effect of density on preferential flow path formation},
  author={Hyman, Jeffrey D},
  journal={Water Resources Research},
  volume={56},
  number={9},
  pages={e2020WR027986},
  year={2020},
  publisher={Wiley Online Library}
}

@article{zhu2021impact,
  title={Impact of fracture geometry and topology on the connectivity and flow properties of stochastic fracture networks},
  author={Zhu, Weiwei and Khirevich, Siarhei and Patzek, Tad W},
  journal={Water Resources Research},
  volume={57},
  number={7},
  pages={e2020WR028652},
  year={2021},
  publisher={Wiley Online Library}
}

@article{moska2021hydraulic,
  title={Hydraulic fracturing in enhanced geothermal systems—field, tectonic and rock mechanics conditions—a review},
  author={Moska, Rafa{\l} and Labus, Krzysztof and Kasza, Piotr},
  journal={Energies},
  volume={14},
  number={18},
  pages={5725},
  year={2021},
  publisher={MDPI}
}

@article{jia2022hydraulic,
  title={Hydraulic stimulation strategies in enhanced geothermal systems (EGS): a review},
  author={Jia, Yunzhong and Tsang, Chin-Fu and Hammar, Axel and Niemi, Auli},
  journal={Geomechanics and Geophysics for Geo-Energy and Geo-Resources},
  volume={8},
  number={6},
  pages={211},
  year={2022},
  publisher={Springer}
}

@article{flores2017effect,
  title={The effect of salinity and gas saturation of a geothermal fluid on the reservoir permeability reduction},
  author={Flores, Ju{\'a}n and Meza, On{\'e}simo and Moya, Sara L and Arag{\'o}n, Alfonso},
  journal={Geof{\'\i}sica internacional},
  volume={56},
  number={4},
  pages={335--343},
  year={2017},
  publisher={Instituto de Geof{\'\i}sica, UNAM}
}

@article{shahidzadeh2010damage,
  title={Damage in porous media due to salt crystallization},
  author={Shahidzadeh-Bonn, Noushine and Desarnaud, Julie and Bertrand, Fran{\c{c}}ois and Chateau, Xavier and Bonn, Daniel},
  journal={Physical Review E—Statistical, Nonlinear, and Soft Matter Physics},
  volume={81},
  number={6},
  pages={066110},
  year={2010},
  publisher={APS}
}

@inproceedings{fridleifsson2008possible,
  title={The possible role and contribution of geothermal energy to the mitigation of climate change},
  author={Fridleifsson, Ingvar B and Bertani, Ruggero and Huenges, Ernst and Lund, John W and Ragnarsson, Arni and Rybach, Ladislaus},
  booktitle={IPCC scoping meeting on renewable energy sources: Proceedings},
  pages={59--80},
  year={2008},
  organization={Intergovernmental Panel on Climate Change}
}

@article{lamur2017permeability,
  title={The permeability of fractured rocks in pressurised volcanic and geothermal systems},
  author={Lamur, A and Kendrick, JE and Eggertsson, GH and Wall, RJ and Ashworth, JD and Lavall{\'e}e, Y},
  journal={Scientific reports},
  volume={7},
  number={1},
  pages={6173},
  year={2017},
  publisher={Nature Publishing Group UK London}
}

@inproceedings{duran2025mixed,
  title={Mixed-Dimensional Approach for Compositional Multiphase Flow in High-Enthalpy Fractured Geothermal Reservoirs},
  author={Duran, Omar and Lipovac, Veljko and Berre, Inga},
  booktitle={Proceedings of the 50th Workshop on Geothermal Reservoir Engineering Stanford University. Stanford Geothermal Workshop. Stanford},
  year={2025}
}

@article{tariq2025fracture,
  title={Fracture Characterization for New Landfill Sites in Crystalline Bedrock: A Case Study from Rogaland, Southwestern Norway},
  author={Tariq, Bilal and French, Helen Kristine and Polteau, St{\'e}phane and Ansch{\"u}tz, Helgard},
  journal={Rock Mechanics and Rock Engineering},
  volume={58},
  number={6},
  pages={7087--7109},
  year={2025},
  publisher={Springer}
}

@article{rahman2002shear,
  title={A shear-dilation-based model for evaluation of hydraulically stimulated naturally fractured reservoirs},
  author={Rahman, MK and Hossain, MM and Rahman, SS},
  journal={International Journal for Numerical and Analytical Methods in Geomechanics},
  volume={26},
  number={5},
  pages={469--497},
  year={2002},
  publisher={Wiley Online Library}
}

@article{tsypkin2004vapour,
  title={Vapour extraction from a water-saturated geothermal reservoir},
  author={Tsypkin, George G and Woods, Andrew W},
  journal={Journal of Fluid Mechanics},
  volume={506},
  pages={315--330},
  year={2004},
  publisher={Cambridge University Press}
}

@article{zhang2023autonomous,
  title={Autonomous fracture flow tunning to enhance efficiency of fractured geothermal systems},
  author={Zhang, Qitao and Taleghani, Arash Dahi},
  journal={Energy},
  volume={281},
  pages={128163},
  year={2023},
  publisher={Elsevier}
}

@article{lei2024thermo,
  title={A thermo-hydro-mechanical simulation on the impact of fracture network connectivity on the production performance of a multi-fracture enhanced geothermal system},
  author={Lei, Zhihong and Zhang, Yulong and Lin, Xingjie and Shi, Yu and Zhang, Yunhui and Zhou, Ling and Shen, Yaping},
  journal={Geothermics},
  volume={122},
  pages={103070},
  year={2024},
  publisher={Elsevier}
}

@article{lu2018global,
  title={A global review of enhanced geothermal system (EGS)},
  author={Lu, Shyi-Min},
  journal={Renewable and Sustainable Energy Reviews},
  volume={81},
  pages={2902--2921},
  year={2018},
  publisher={Elsevier}
}

@article{shukla2022climate,
  title={Climate change 2022: Mitigation of climate change},
  author={Shukla, Priyadarshi R and Skea, Jim and Slade, Raphael and Al Khourdajie, Alaa and van Diemen, Ren{\'e}e and McCollum, David and Pathak, Minal and Some, Shreya and Vyas, Purvi and Fradera, Roger and others},
  journal={Contribution of working group III to the sixth assessment report of the Intergovernmental Panel on Climate Change},
  volume={10},
  pages={9781009157926},
  year={2022},
  publisher={Cambridge University Press Cambridge, UK}
}

@article{fridleifsson2001geothermal,
  title={Geothermal energy for the benefit of the people},
  author={Fridleifsson, Ingvar B},
  journal={Renewable and sustainable energy reviews},
  volume={5},
  number={3},
  pages={299--312},
  year={2001},
  publisher={Elsevier}
}

@article{wang2021modeling,
  title={Modeling of multiphase mass and heat transfer in fractured high-enthalpy geothermal systems with advanced discrete fracture methodology},
  author={Wang, Yang and de Hoop, Stephan and Voskov, Denis and Bruhn, David and Bertotti, Giovanni},
  journal={Advances in Water Resources},
  volume={154},
  pages={103985},
  year={2021},
  publisher={Elsevier}
}

@book {nocedal2006Numerical,

	title = {Numerical Optimization},

	publisher = {Springer},

	author = {Nocedal, J. and Wright, S. J.},

	edition = {Second},

	year = 2006

}

@inproceedings{oguntola2020robust,
	title={On the Robust Value Quantification of Polymer {EOR} Injection Strategies for Better Decision Making},
	author={Oguntola, M. and Lorentzen, R.},
	booktitle={ECMOR XVII},
	volume={2020},
	pages={1--25},
	year={2020},
	organization={European Association of Geoscientists \& Engineers}
}

@Article{aavatsmark2002introduction,
	title   = {An introduction to multipoint flux approximations for quadrilateral grids},
	author  = {Aavatsmark, I.},
	journal = {Computational Geosciences},
	volume  = {6},
	pages   = {405--432},
	year    = {2002}
}

@Article{aghili2020hybrid,
	title   = {A hybrid-dimensional compositional two-phase flow model in fractured porous media with phase transitions and Fickian diffusion},
	author  = {Aghili, J. and Dreuzy, J. de and Masson, R. and Trenty, L.},
	journal = {Journal of Computational Physics},
	volume  = {441},
	pages   = {110452},
	year    = {2020},
	doi     = {10.1016/j.jcp.2021.110452}
}

@Article{alpak2018variable,
	title   = {A variable-switching method for mass-variable-based reservoir simulators},
	author  = {Alpak, F. O. and Vink, J. C.},
	journal = {Spe Journal},
	volume  = {23},
	number  = {5},
	pages   = {1469--1495},
	year    = {2018},
	doi     = {10.2118/182606-pa}
}

@Article{beaude2019non,
	title   = {Non-isothermal compositional liquid gas Darcy flow: Formulation, soil-atmosphere boundary condition and application to high-energy geothermal simulations},
	author  = {Beaude, L. and Brenner, K. and Lopez, S. and Masson, R. and Smaï, F.},
	journal = {Computational Geosciences},
	volume  = {23},
	number  = {3},
	pages   = {1--28},
	year    = {2019},
	doi     = {10.1007/s10596-018-9794-9}
}

@article{schneider2020coupling,
  title={Coupling staggered-grid and MPFA finite volume methods for free flow/porous-medium flow problems},
  author={Schneider, Martin and Weishaupt, Kilian and Gl{\"a}ser, Dennis and Boon, Wietse M and Helmig, Rainer},
  journal={Journal of Computational Physics},
  volume={401},
  pages={109012},
  year={2020},
  publisher={Elsevier}
}

@Book{chen2006computational,
	title     = {Computational methods for multiphase flows in porous media},
	author    = {Chen, Z. and Huan, G. and Ma, Y.},
	year      = {2006},
	publisher = {SIAM},
	address   = {Philadelphia}
}

@article{wang2022high,
  title={High-enthalpy geothermal simulation with continuous localization in physics},
  author={Wang, Yang and Voskov, Denis},
  journal={Mathematics},
  volume={10},
  number={22},
  pages={4328},
  year={2022},
  publisher={MDPI}
}

@article{voskov2017operator,
  title={Operator-based linearization approach for modeling of multiphase multi-component flow in porous media},
  author={Voskov, Denis V},
  journal={Journal of Computational Physics},
  volume={337},
  pages={275--288},
  year={2017},
  publisher={Elsevier}
}

@Article{class2002numerical,
	title   = {Numerical simulation of non-isothermal multiphase multicomponent processes in porous media. 2. Applications for the injection of steam and air},
	author  = {Class, H. and Helmig, R.},
	journal = {Advances in Water Resources},
	volume  = {25},
	number  = {5},
	pages   = {551--564},
	year    = {2002}
}

@Article{driesner2007system,
	title   = {The system H2O--NaCl. Part I: Correlation formulae for phase relations in temperature-pressure-composition space from 0 to 1000 C, 0 to 5000 bar, and 0 to 1 XNaCl},
	author  = {Driesner, T. and Christoph, H.},
	journal = {Geochimica et Cosmochimica Acta},
	volume  = {71},
	number  = {20},
	pages   = {4880--4901},
	year    = {2007}
}

@article{lipovac2025persistent,
  title={Persistent-variable thermal compositional simulation of multiphase flow with phase separation in porous media},
  author={Lipovac, Veljko and Duran, Omar and Keilegavlen, Eirik and Berre, Inga},
  journal={arXiv preprint arXiv:2512.04205},
  year={2025}
}

@Article{falko2021brine,
	title   = {Brine formation and mobilization in submarine hydrothermal systems: Insights from a novel multiphase hydrothermal flow model in the system H2O--NaCl},
	author  = {Falko, V. and Jörg, H. and Lars, R.},
	journal = {Transport in Porous Media},
	volume  = {136},
	pages   = {65--102},
	year    = {2021}
}

@Article{geiger2006multiphase1,
	title   = {Multiphase thermohaline convection in the earth's crust: I. A new finite element—Finite volume solution technique combined with a new equation of state for NaCl-H2O},
	author  = {Geiger, S. and Driesner, T. and Heinrich, C. A. and Matthäi, S.},
	journal = {Transport in Porous Media},
	volume  = {63},
	number  = {3},
	pages   = {399--434},
	year    = {2006},
	doi     = {10.1007/s11242-005-0108-z}
}

@Article{geiger2006multiphase2,
	title   = {Multiphase thermohaline convection in the earth's crust: II. Benchmarking and application of a finite element—Finite volume solution technique with a NaCl-H2O equation of state},
	author  = {Geiger, S. and Driesner, T. and Heinrich, C. A. and Matthäi, S.},
	journal = {Transport in Porous Media},
	volume  = {63},
	number  = {3},
	pages   = {435--461},
	year    = {2006},
	doi     = {10.1007/s11242-005-0109-y}
}

@inproceedings{oguntola2025unified,
  title     = {A Unified Compositional Flow Model for Simulating Multiphase High-Enthalpy Geothermal Reservoirs},
  author={Oguntola, Micheal Babatunde and Duran, Omar and Lipovac, Veljko and Keilegavlen, Eirik and Berre, Inga},
  booktitle = {Proceedings of the 50th Workshop on Geothermal Reservoir Engineering},
  year      = {2025},
  address   = {Stanford, California},
  month     = {February 10--12},
  organization = {Stanford University}
}

@Article{keilegavlen2021porepy,
	title   = {PorePy: An open-source software for simulation of multiphysics processes in fractured porous media},
	author  = {Keilegavlen, E. and Berge, R. L. and Fumagalli, A. and Starnoni, M. and Stefansson, I. and Varela, J. and Berre, I.},
	journal = {Computational Geosciences},
	volume  = {25},
	pages   = {243--265},
	year    = {2021},
	doi     = {10.1007/s10596-020-10002-5}
}

@Misc{khait2019delft,
	title        = {Delft Advanced Research Terra Simulator (DARTS): General purpose reservoir simulator with operator-based linearization},
	author       = {Khait, M.},
	year         = {2019},
	howpublished = {\url{https://doi.org/10.4233/uuid:5f0f9b80-a7d6-488d-9bd2-d68b9d7b4b87}}
}

@Article{kipp2008guide,
	title   = {Guide to the revised ground-water flow and heat transport simulator: HYDROTHERM - Version 3},
	author  = {Kipp, K. L. and Hsieh, P. A. and Charlton, S. R. and Charlton, S. R.},
	journal = {Techniques and Methods},
	year    = {2008},
	doi     = {10.3133/tm6a25}
}

@Article{lauser2011new,
	title   = {A new approach for phase transitions in miscible multi-phase flow in porous media},
	author  = {Lauser, A. and Hager, C. and Helmig, R. and Wohlmuth, B.},
	journal = {Advances in Water Resources},
	volume  = {34},
	number  = {8},
	pages   = {957--966},
	year    = {2011},
	doi     = {10.1016/j.advwatres.2011.04.021}
}

@TechReport{pruess2003tough,
	title       = {The TOUGH codes—A family of simulation tools for multiphase flow and transport processes in permeable media},
	author      = {Pruess, K.},
	institution = {Lawrence Berkeley National Laboratory},
	year        = {2003},
	doi         = {10.2136/vzj2004.0738}
}

@Article{quiroz2024multi,
	title   = {Multi-segmented non-isothermal compositional liquid gas well model for geothermal processes},
	author  = {Quiroz, D. C. and Jeannin, L. and Lopez, S. and Masson, R.},
	journal = {arXiv.Org},
	year    = {2024},
	doi     = {10.48550/arxiv.2401.02406}
}

@Article{rajabi2023dynamical,
	title   = {Dynamical modeling of a geothermal system to predict hot spring behavior},
	author  = {Rajabi, M. M. and Chen, M.},
	journal = {Modeling Earth Systems and Environment},
	volume  = {9},
	number  = {3},
	pages   = {3085--3093},
	year    = {2023},
	doi     = {10.1007/s40808-023-01696-4}
}

@Article{voskov2012comparison,
	title   = {Comparison of nonlinear formulations for two-phase multi-component EoS based simulation},
	author  = {Voskov, D. and Tchelepi, H.},
	journal = {Journal of Petroleum Science and Engineering},
	volume  = {82},
	pages   = {101--111},
	year    = {2012}
}

@Article{weis2014hydrothermal,
	title   = {Hydrothermal, multiphase convection of H2O‐NaCl fluids from ambient to magmatic temperatures: A new numerical scheme and benchmarks for code comparison},
	author  = {Weis, P. and Driesner, T. and Coumou, D. and Geiger, S.},
	journal = {Geofluids},
	volume  = {14},
	number  = {3},
	pages   = {347--371},
	year    = {2014},
	doi     = {10.1111/gfl.12080}
}

@article{stefansson2021fully,
  title={A fully coupled numerical model of thermo-hydro-mechanical processes and fracture contact mechanics in porous media},
  author={Stefansson, Ivar and Berre, Inga and Keilegavlen, Eirik},
  journal={Computer Methods in Applied Mechanics and Engineering},
  volume={386},
  pages={114122},
  year={2021},
  publisher={Elsevier}
}

\end{document}